\documentclass[preprint,12pt]{elsarticle}
\usepackage[margin=.8in]{geometry}

\usepackage{amsmath,amssymb,float,color,amsthm,caption}

\usepackage{graphicx} 
\usepackage{booktabs} 
\usepackage{array} 
\usepackage{verbatim} 
\usepackage{psfrag}
\usepackage{braket}

 \graphicspath{{figure/}}

\usepackage{marginnote}

\newcommand{\lbl}[1]{\label{#1}}

\def\grad{{\nabla}}
\def\div{{\nabla\cdot}}
\def\half{{\textstyle \frac{1}{2}}}
\def\cM{\mathcal{M}}
\def\cW{\mathcal{W}}
\def\cG{\mathcal{G}}
\def\veps{{\varepsilon}}

\def\x{{\bf x}}

\def\MB{\overline{\cM}}  
\def\Mb{M_2}  
\def\Mh{M_1}  
\def\Fe{F_{\mathrm{ex}}}
\def\Fi{F_{\mathrm{im}}}
\def\tf{T}

\usepackage{scalerel,stackengine}
\stackMath
\newcommand\reallywidehat[1]{%
\savestack{\tmpbox}{\stretchto{%
  \scaleto{%
    \scalerel*[\widthof{\ensuremath{#1}}]{\kern.1pt\mathchar"0362\kern.1pt}%
    {\rule{0ex}{\textheight}}
  }{\textheight}%
}{2.4ex}}%
\stackon[-6.9pt]{#1}{\tmpbox}%
}
\usepackage{graphicx}
\usepackage{epstopdf}

\journal{Computational Materials Science}

\begin{document}

\begin{frontmatter}



\title{IMEX methods for thin-film equations and Cahn-Hilliard equations with variable mobility}


\author[label1]{Saulo Orizaga\corref{cor1}}
\author[label2]{Thomas Witelski}
\tnotetext[cor1]{corresponding author: saulo.orizaga@nmt.edu}

\affiliation[label1]{organization={ Department of Mathematics, New Mexico Institute of Mining and Technology},
            addressline={801 Leroy Place}, 
            city={Socorro},
            postcode={87801}, 
            state={NM},
            country={USA}}
            
\affiliation[label2]{organization={ Department of Mathematics, Duke University},
            addressline={120 Science Dr}, 
            city={Durham},
            postcode={27708-0320}, 
            state={NC},
            country={USA}}

\begin{abstract}

%
We explore a class of splitting schemes employing implicit-explicit (IMEX) time-stepping to achieve accurate and energy-stable solutions for thin-film equations and Cahn-Hilliard models with variable mobility. This splitting method incorporates a linear, constant coefficient implicit step, facilitating efficient computational implementation. We investigate the influence of stabilizing splitting parameters on the numerical solution computationally, considering various initial conditions. Furthermore, we generate energy-stability plots for the proposed methods, examining different choices of splitting parameter values and timestep sizes. These methods enhance the accuracy of the original bi-harmonic-modified (BHM) approach, while preserving its energy-decreasing property and achieving second-order accuracy. We present numerical experiments to illustrate the performance of the proposed methods.\\
\end{abstract}

\begin{graphicalabstract}
\includegraphics[width=\textwidth]{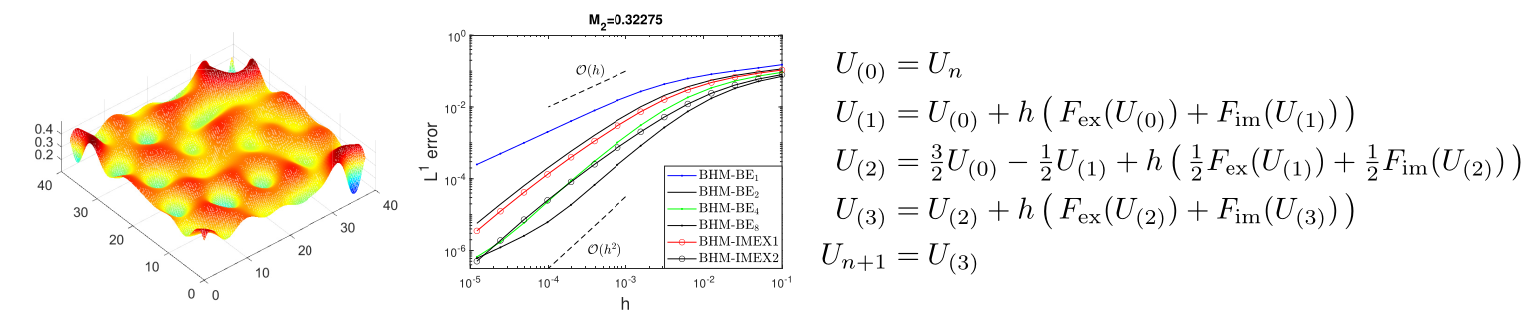}

\end{graphicalabstract}

\begin{highlights}
\item Splitting methods for Cahn-Hilliard (CH) equations with variable mobility.
\item Efficient and accurate numerical schemes for thin film (TF) equations.
\item Testing of energy-stability for schemes for TF and CH equations.
\item  Study of second and fourth order splitting parameters for TF equations.
\end{highlights}


\end{frontmatter}


\section{Introduction}
%
Many problems in materials science, mathematical biology, and other areas can be described in terms of phase-field models. These models use higher-order parabolic partial differential equations to describe the evolution of a phase function $u({\bf x},t)$. In this paper, we will focus on numerical methods for models of the form
\begin{equation}
{\partial u \over \partial t} = \nabla \cdot
(\cM(u)\nabla [\cW'(u) -\Delta u]).
    \lbl{CHVMeqn}
\end{equation}
The integral of $u$ is conserved thanks to the divergence form of this evolution equation and 
the flux is the product of a non-negative mobility function $\cM(u)\ge 0$ and the gradient of a 
chemical potential $w$. This potential is the variational derivative of an energy functional,
$w\equiv \frac{\delta \mathcal{E}}{\delta u} = \cW'(u)-\Delta u$ with 
\begin{equation}
\mathcal{E}(u) = \int_{\Omega} \cW(u)+ \frac{1}{2} |\grad u|^2  \, dx \,.
\lbl{defenergy}
\end{equation}
It follows that
the dynamics of \eqref{CHVMeqn} is a gradient flow in a weighted $H^{-1}$ norm for the monotone decreasing energy functional 
\begin{equation}
\frac{d\mathcal{E}}{dt}  = - \int_{\Omega}  \cM(u) \left|  \grad w \right|^2 \, dx \leq 0\,.
\lbl{dEdt}
\end{equation}
Setting a constant mobility $\cM=\epsilon^2$ and potential $\cW(u)=(u^4-2u^2)/(4\epsilon^2)$ yields the classic Cahn-Hilliard (CH) equation,
\begin{equation}
{\partial u \over \partial t} = \Delta( u^3-u -\epsilon^2 \Delta u),
\lbl{CHclassic}
\end{equation}
where $\epsilon>0$ defines the lengthscale for the width of transition
layers between equilibrium phases $u= \pm 1$.
This equation was originally formulated to model phase separation in a binary alloy mixture \cite{cahn1958free} with the range $-1\le
u\le 1$ defining physically meaningful mixtures. Since then it has been applied to many other problems including block co-polymers, image processing, and biological systems  \cite{miranville2019,ch-app1, ch-app2, ch-app3,garcke2017well}.

\par

There are extensive bodies of research on the
analysis \cite{miranville2019,ANC1984,miranville2017cahn} and efficient
numerical computations for this equation \cite{glasnerOrizaga}.
Generic initial conditions starting near the unstable state $u=0$ rapidly
separate (sometimes called spinodal decomposition) into complicated arrays
of nearly pure ($u=\pm 1$) local states followed by slow dynamics of the
phase boundaries leading to merging and re-arrangements of the pure states
(called coarsening or Ostwald ripening). If the initial data satisfies $|u_0|<1$, then
having the solution remain in the range $-1\le u \le 1$ for all times
has been called the ``positivity property'' \cite{chen2019pos}.  This is known
to hold in one dimension, but not more generally due to the curvature of phase interfaces (sometimes called the Gibbs-Thomson effect), see for example \cite{DaiDu1,DaiDu2}.
\par 
The Cahn-Hilliard equation with concentration-dependent mobility $\cM(u)$ is
known to yield more detailed models for many physical applications in terms
of better capturing the timescales of the dynamics. The particular forms of $\cM(u)$ that are of most interest
exhibit degeneracy, reducing the flux as one or both of the pure states are
approached, $\cM(u)\to 0$ as $u\to \pm 1$. Such degeneracy makes mathematical
analysis of these models more challenging for obtaining results on existence
and uniqueness of solutions \cite{ElliottGarcke}. Indeed, in this case, it
is not guaranteed that solutions remain in the physical range, $|u|\le 1$,
\cite{DaiDu2} and consequently the design of numerical methods must address
how to handle such possible behaviors \cite{DaiDu1,HectorC2013}.  
Very importantly, if positivity is not preserved then forms for $\cM(u)$
that have negative mobilities for $u$ outside $|u|\le 1$ cannot be used
since they would allow the problem to become illposed.

\par
Another important class of applications fitting within the framework of \eqref{CHVMeqn}
 are models of viscous thin film fluid dynamics \cite{Thiele2018}.
Derived from lubrication theory in the asymptotic limit of a low Reynolds
number and small aspect ratio, thin film (TF) equations describe the evolution of the
thickness, $u({\bf x},t)\ge 0$, of fluid layers coating solid
substrates \cite{oron,craster2009dynamics}.  In scaled
non-dimensional form, the model is a nonlinear fourth-order parabolic
diffusion equation for $u$, of the form
\begin{equation}
   \lbl{tfe}
u_t =\div\left(u^3 \grad [\Pi(u) -\Delta u]\right)\,. \\
\end{equation}
Here the role of the chemical potential is played
by the hydrostatic pressure,
$p=\Pi(u)-\Delta u$, giving pressure contributions from surface tension and
conservative forces describing physiochemical material properties.  The $\Pi(u)$ function, called a disjoining pressure, is the derivative of a molecular potential energy, as in $\Pi(u)=\cW'(u)$.
The mobility in \eqref{tfe}, $\cM(u)=u^3$, is widely used for classic models of fluid mechanics, but  other forms of degenerate mobilities have been considered such as $\cM(u)=u^n$ for
$n>0$ \cite{bertozzi1998}. 
One simple form for $\Pi(u)$ modeling fluid spreading on a hydrophobic solid is
\begin{equation}
\Pi(u)=\frac{\veps^2}{u^3}-\frac{\veps^3}{u^4}
\lbl{eqPi}
\end{equation}
where  $\veps>0$  sets a minimal
film thickness determined by materials properties
\cite{Bertozzi2001}. This gives a double-well
potential comparable to CH models generating phase separation between 
nearly-uniform ``precursor layers'' with $u=O(\veps)\ll 1$ for
$\veps\to 0$ and finite mass fluid droplets \cite{glasnertom}.
In \cite{Bertozzi2001} it was proved for the one-dimensional case that the
singular nature of \eqref{eqPi} for $u\to 0$ is sufficient to ensure that
$u$ maintains positivity ($u> 0$) for all times in \eqref{tfe}, avoiding the
difficulties \cite{HectorC2013} for the degenerate-mobility CH equation when
the potential is non-singular. 
The stages of dynamics for thin film equations share most of the richness and complexity of the dynamics
of other Cahn-Hilliard problems \cite{witelski2020}. While CH models in
materials science and other settings are applied to two- and
three-dimensional problems, since TF equations describe fluid coatings,
their applications are limited to two dimensions. 
A key difference that can make computing TF problem more challenging is that while CH solutions nominally span a finite range of values, $|u|\lesssim 1$, TF solutions have no a priori upper bound and hence the influence of the variable mobility can play a stronger role. 
The effective $\cW(u)$ potential for TF problems has a degenerate second well at $u\to\infty$.
For both CH and TF problems, computations are generally used to describe statistical properties on large domains evolving over long-time simulations. This has motivated studies to advance
efficient computational methods for these problems
\cite{LAM2019,kondic2002,adi2003}. 

\par
Numerical methods with explicit time discretization are impractical for these equations due to time-stepping restrictions associated with numerical stability.  This has motivated the development of a variety of
implicit and semi-implicit numerical methods
\cite{ascher1995implicit,chen1998applications,eyre1998unconditionally,he2007large,gomez2011provably,cheng2008efficient,wise2009energy}.  
A simple and elegant semi-implicit approach was formulated by the highly influential work of Eyre \cite{eyre1998unconditionally}. For the general phase-field model with variable mobility \eqref{CHVMeqn}, applying Eyre's splitting yields
\begin{equation}
{\partial u\over \partial t} =\div\left( \cM(u) \grad \left[\frac{\delta \mathcal{E_+}}{\delta u} + \frac{\delta \mathcal{E_-}}{\delta u}\right]\right) 
\lbl{eyre_original}
\end{equation}
where the energy is separated into convex and concave parts $ \mathcal{E}=\mathcal{E_+}+\mathcal{E_-} $ and the convex and concave parts of the energy are to be computed
implicitly and explicitly respectively. If the mobility function is taken to be constant $(\cM(u)=1)$  then \eqref{eyre_original} simplifies to
\begin{equation}
{\partial u\over \partial t} =\Delta^2 \left[\frac{\delta \mathcal{E_+}}{\delta u} + \frac{\delta \mathcal{E_-}}{\delta u}\right].
\lbl{eyre_original2}
\end{equation}
It is often the case that the convex part of the energy $  \Delta^2  \left[\frac{\delta \mathcal{E_+}}{\delta u} \right] $ is chosen to be linear in $u$ so this term  
 can be computed implicitly very efficiently without the use of non-linear solvers  associated with
Newton iterations.
For the classic Cahn-Hilliard equation \eqref{CHclassic}, the energy is separated into convex and concave parts which implies $\cW(u)=\cW_+(u)+\cW_-(u)$ with $\cW_+''(u)>0$ and $\cW_-''(u)<0$,
\begin{equation}
{\partial u\over \partial t} = \underbrace{\Delta \cW_+'(u) -\epsilon^2 \Delta^2 u}_{F_{\mathrm{contract}}(u)} \, 
+\, \underbrace{\Delta \cW_-'(u)}_{F_{\mathrm{expans}}(u)}
\lbl{eyre}
\end{equation}
yielding the indicated contractive and expansive operators. 
In numerical implementations $F_{\mathrm{contract}}$ is treated implicitly and $F_{\mathrm{expans}}$ explicitly.
Convexity-splitting has been successfully implemented in
many different contexts \cite{glasner2006grain,greer2006fourth,wise2009energy,cristini2009nonlinear}.
Various improvements and extensions of the method have also been made 
\cite{boyer2011numerical,zhang2013adaptive,han2015second}.
For the class of problems given by equation (\ref{CHclassic}) and the case of
constant mobility, a series of schemes based on the convexity
splitting idea were proposed by Glasner and Orizaga \cite{glasnerOrizaga}.

\par
On the other hand, if  $\cM(u)$ is variable then nonlinear terms in 
$ \div\left( \cM(u) \grad \left[\frac{\delta \mathcal{E_+}}{\delta u}\right] \right) $ 
seem unavoidable.
Unfortunately, the methods proposed in
\cite{eyre1998unconditionally,glasnerOrizaga} are not applicable to the 
above-mentioned equations due to the concentration-dependent mobility $\cM(u)$. 
This
gives the motivation for proposing a new framework that consists of a
first-order energy-stable splitting approach coupled with an
implicit-explicit (IMEX) time-stepping discretization to obtain accurate and
energy-stable solutions for equations of the form (\ref{CHVMeqn}).

\par
We will refer to \eqref{CHVMeqn} as the Cahn-Hilliard equation with variable mobility (CHVM) 
and use it as the framework for considering numerical methods for both Cahn-Hilliard and thin-film type problems.
For some problems it will be convenient to re-expand this equation 
to separate the second- and fourth-order spatial operators, which yields
\begin{equation}
\lbl{nlfo2}
{\partial u \over \partial t} = \underbrace{\Delta \cG(u) -\nabla\cdot (\cM(u) \grad \Delta u\vphantom{\int_a^b}
)}_{F(u)}, 
\qquad
\mbox{where} 
\qquad
\cG(u)=\int \cM(u)\cW''(u)\,du\,.
\end{equation}
This form allows us to consider approaches splitting and stabilization that deal separately with the nonlinearities of the second- and fourth-order terms.
We will write models \eqref{CHVMeqn} and \eqref{nlfo2} in the general form
\begin{equation}
{\partial u\over \partial t} = F(u) \qquad \implies \qquad 
{\partial u\over \partial t} = \Fi(u)+\Fe(u)
\lbl{odeF}
\end{equation}
with splitting of the spatial operator into a linear stabilizing part that will be treated implicitly
\begin{equation}
\Fi(u) = -\sum_{i=0}^2 M_i (-\Delta)^i u
\lbl{FimEqn}
\end{equation}
with constant coefficients $M_i$
and the remainder that will be handled explicitly, $\Fe(u)=F(u)-\Fi(u)$.

\par
In the study of time-dependent problems, splitting methods have been shown to provide efficient computational results.
Different splittings methods have been proposed, for the class of equations of the form (\ref{CHVMeqn}), in order for numerical implementations to produce
accurate and stable solutions. Some of the early work related to splitting
methods can be traced back to early 1970's in the work of Dupont and Douglas
\cite{douglasdupont}.  More recently, two splittings methods can be found in
the work by Barrett
and Blowey (1999) \cite{barrettblowey} and Bertozzi et al.\ (2011)
\cite{bertozzi}. These methods have similar ideas, but the main distinction
is that \cite{bertozzi} proposes a scheme that is linear
with respect to its implicit terms
while in \cite{barrettblowey} they consider a nonlinear implicit scheme that
likewise builds on the constant mobility case for the highest-order terms.
In \cite{eggers2014} Duchemin and Eggers used linear stability analysis to show that forward Euler schemes can always be stabilized with a splitting having a single-term $\Fi(u)$ with an appropriately selected coefficient. Their results are clear and broadly applicable and provide necessary conditions to ensure stability. In considering the more general form \eqref{FimEqn} our focus is primarily on how the $M_i$ coefficients impact the accuracy of the scheme.

\par
For fourth-order phase-field models, choices for the splitting parameters $\Mb, \Mh$ in \eqref{FimEqn} have been discussed in many previous papers 
\cite{smereka,hongkai2003,barrettblowey,du2015,shenSAV2018,jiang2012,li2016,yanwise2018,zhutikare1999} but remain a very important issue with many unresolved questions. Values have often been selected to ensure convexity splitting (as in \eqref{eyre}) or based on proofs for energy stability, so, like \eqref{dEdt}, the discretized dynamics have monotone decreasing energy, $\mathcal{E}(U_{n+1})\le \mathcal{E}(U_n)$. However, some of these papers also comment on the parameter values not being definitive \cite{jiang2012,yanwise2018}. We will re-examine and interpret the influences of $\Mb, \Mh$ on energy-stability and accuracy of the numerical methods. We will consider different cases for the splitting parameter (fixed constant) as it is typically done for CH-type equations in which solutions are bounded on $[-1,1]$ \cite{glasnerOrizaga}. For TF equations the solution has a lower bound $\epsilon$ (thin layer of fluid), but the solution maximum can exceed the value of $u \approx 1$. For this reason, a fixed value of a splitting parameter may not be suitable for TF equations. 
We will examine results when these splitting parameters are fixed constants and also when they are allowed to evolve dynamically in relation to properties of the solution, as considered in \cite{bertozzi}.

\par
The structure of the paper is as follows: Section 2 discusses the biharmonic modified
method applied to time-dependent problems and introduces the four main numerical methods considered in this paper BHM-BE$_J$,
BHM-CN$_J$, BHM-IMEX1 and BHM-IMEX2. In Section 3, we present results from numerical simulations to illustrate
issues for the order of accuracy of the schemes and dependence on the splitting parameter as well as considerations of qualitative properties of
numerical solutions regarding energy-stability. Section 4, we provide illustrations of computed dynamics for both the thin-film equation and the Cahn-Hilliard
equation with variable mobility for long enough running times to demonstrate the performance of the methods. 
%



\section{Methods}
\subsection{Biharmonic splitting and extensions}

\par
Applied to \eqref{nlfo2}, the splitting method proposed in \cite{bertozzi} can be constructed
 in the following way. The mobility
coefficient function in the fourth-order operator can be written as
$\cM(u)=\Mb+(\cM(u)-\Mb)$ to give the form
\begin{equation}
\lbl{nlfo3}
{\partial u\over \partial t} = 
\underbrace{-\Mb \Delta^2 u}_{\Fi(u)}\;+\; \underbrace{\Delta \cG(u)-\div
[\left(\cM(u)-\Mb\right) \grad \Delta u ]}_{\Fe(u)}.
\end{equation}
This fits the form (\ref{odeF}, \ref{FimEqn}) setting $\Mh=M_0=0$.
We denote time-discretized approximation of the solution as
$u({\bf x},t_n)\approx U_n$ where the discrete times will be
expressed with respect to the local timestep, $h$, by $t_{n+1}=t_n+h$.
Using the backward Euler difference for the time derivative and the $\Mb$
term in \eqref{nlfo3} and treating $\Fe(u)$ explicitly yields the first-order
implicit-explicit semi-discrete method,
\begin{equation}
\lbl{schemeBH}
\frac{U_{n+1}-U_{n}}{h}=\Fi(U_{n+1}) + \Fe(U_n),
\qquad \mbox{(BHM)}
\end{equation}
called the biharmonic modified scheme in \cite{bertozzi}. If no stabilization is applied, $\Fi\equiv 0$ and
$\Fe=F(u)$, and \eqref{schemeBH} reduces to the explicit forward Euler scheme for \eqref{odeF}.
Our focus is on choices for different forms for the splitting of implicit and explicit
terms in the time-stepping schemes that can be applied with general spatial
discretizations such as spectral methods, finite elements, or finite
differences.

\par
We will use Fourier pseudo-spectral methods \cite{trefethen2000spectral} for our computations. 
Fourier methods have been widely applied to Cahn-Hilliard models in many studies \cite{shengdu2010,zhutikare1999,li2016,li2017,xu2019}. We illustrate the spatial discretization for \eqref{schemeBH} (and the application to other schemes we consider follows similarly). The two-dimensional discrete Fourier transform of $U$ can be written as
$$
U \approx \sum_{k_x= -N/2}^{N/2-1} \sum_{k_y= -N/2}^{N/2-1}   \widehat{U}(k_x,k_y,t) \exp \Big[i\left( k_x x +k_y y   \right)     \Big]. 
$$
Consequently, the Fourier transforms of the harmonic and biharmonic operators can be expressed as
$\widehat{\Delta U} =-k^2 \widehat U$ and $\widehat{\Delta^2 U} =k^4 \widehat U$ where $k^2=k_x^2+k_y^2$. Applied to \eqref{schemeBH}, this yields
\begin{equation}
\widehat U_{n+1}=\frac{\widehat U_{n} +h\reallywidehat{\Fe(U_n)} }{1+h \Mb k^4},
\lbl{spectral}
\end{equation}
giving an efficient means for evaluating $\widehat{U}_{n+1}$ and then the inverse transform is applied to obtain $U_{n+1}$. Following the approach in \cite{eggers2014}, this form can be used to analyze conditions for linear stability by re-writing $\reallywidehat{\Fe(U_n)}=\reallywidehat{F(U_n)}-\reallywidehat{\Fi(U_n)}$ yielding
\begin{equation}
\widehat U_{n+1}=\widehat U_{n} +\frac{h\reallywidehat{F(U_n)} }{1+h \Mb k^4}.
\end{equation}
Then the principal part (highest order terms) of $F$ can be used to estimate $\reallywidehat{F(U_n+e_n)}\approx\reallywidehat{F(U_n)} -C_0 k^4 \widehat{e}_n$ and the amplification factor for the growth of errors to $U_n$ is given by $\sigma=1-hC_0 k^4/(1+h M_2 k^4)$. Unconditional stability can then be achieved if $M_2> C_0/2$.

\par
The local truncation error for \eqref{schemeBH} can be derived by Taylor-expanding $U(t)$ about $t_n$ for $h\to 0$ to yield
\begin{equation}
\tau_n = \left[\half {\delta F_n\over \delta u} F(U_n) + M_2 \Delta^2 F(U_n)\right] h + O(h^2)
\lbl{TruncErr}
\end{equation}
where $\delta_u F_n$ is the functional derivative of $F(u)$ evaluated at $U_n$. As is expected, $\tau_n$ is formally first-order accurate, and we can see the magnitude of the error depends on how the full operator compares with the biharmonic stabilizing term. 
As described in \cite{eggers2014}, while linearized stability can always be achieved in this scheme, considerations of the accuracy of the computed solution (i.e. the size of the coefficient in $\tau\sim C_1 h$) may become the limiting factor for determining time-steps.

\par
A generalization of Bertozzi's backward Euler 
biharmonic splitting \eqref{schemeBH} can be
written as 
\begin{equation}
    {U_{(j)} -  U_n \over h} -\Fi(U_{(j)}) = \Fe(U_{(j-1)})
\hskip1in \mbox{(BHM-BE$_J$)}
\lbl{BHMBE}
\end{equation}
for $j=1,2,\cdots, J$,
where the initial iterate at each timestep can be $U_{(0)}=U_n$
(or some multistep extrapolation when available, like $U_{(0)}=2U_n-U_{n-1}$
\cite{glasnerOrizaga}). Without iteration, this reduces to
Bertozzi's scheme \eqref{schemeBH},
with $U_{(0)}=U_n$ and generating $U_{n+1}=U_{(1)}$.
By iterating \eqref{BHMBE} at each time-step, more accurate
solutions can be obtained. To see that the iterates
for $j=1,2,3,\cdots,J$ approach the backward Euler
solution, $U_{(j)}\to U_{n+1}$ if they converge, write $U_{(j)}=U_{n+1} + e_{(j)}$. Using fixed-point analysis, the sequence will converge, with $|e_{(j)}|<|e_{(j-1)}|$, if $||h{({\bf I} +hM_2\Delta^2)}^{-1} \delta_u \Fe(U_{n+1})||<1$.

\par
Similarly, we can formulate a second-order Crank-Nicolson (or trapezoidal) type scheme as
\begin{equation}
    {U_{(j)} -  U_n \over h} -\half \Fi( U_{(j)}) =
\half  \Fe(U_{(j-1)})+ \half F(U_n)
\qquad \mbox{(BHM-CN$_J$)}
\lbl{BHMCN}
\end{equation}
for $j=1,2,\cdots, J$. Without iteration, CN$_1$ is first-order accurate and reduces to \eqref{schemeBH} with $M_{2,CN}=2M_{2,BE}$. We will see that with iterations, the solutions approach second-order accuracy 
with the error coefficient approaching the lower bound set by the results
from the true Crank-Nicolson nonlinear implicit scheme.

\par
More generally, both second-order and fourth-order splitting can be considered with the above time-stepping schemes,
\begin{equation}
\lbl{nlfo4}
{\partial u\over \partial t} = 
\underbrace{-\Mb \Delta^2 u+\Mh\Delta u}_{\Fi(u)}\;+\; \underbrace{
\Delta(\cG(u)-\Mh u)-\div
[\left(\cM(u)-\Mb\right) \grad \Delta u ]}_{\Fe(u)},  
\end{equation}
Here we will call $\Mb$ the fourth-order splitting or stabilization parameter, and $\Mh$ is the second-order parameter. For $\Mh=0$ this reverts to \eqref{nlfo3}, while retaining $\Mh$ allows for a linear-type convexity splitting of the potential. Since the definition of $G(u)$ combines the potential with the mobility, we note that \eqref{nlfo4} can also be re-written in terms of the linear splitting $\cW'=W_1 u + (\cW'(u)-W_1 u)$ yielding 
$\Fe(u) = \nabla\cdot (\cM(u)\nabla (\cW'(u)-W_1 u))
-\div[\left(\cM(u)-\Mb\right) \grad( \Delta u-W_1 u) ]
$
which is equivalent to \eqref{nlfo4} if $\Mh=\Mb W_1$. For constant mobility, with $\cM(u)\equiv \Mb$, \eqref{nlfo4} reduces to splitting methods applied only to the second-order operators, as considered in \cite{glasnerOrizaga} and other papers.  The task of computing the nonlinear part of energy implicitly in \eqref{eyre_original} is mitigated with the bi-harmonic splitting. 
Equation \eqref{nlfo4} is a modified version of the original Eyre's splitting
\eqref{eyre_original} and for this reason unconditional energy stability may not be guaranteed. 


\subsection{IMEX Methods}
Implicit-explicit (IMEX) schemes can be seen as extensions and improvements over the original Runge-Kutta schemes. 
For our purposes, we can informally understand IMEX schemes as integrating tools that apply to model equations after a splitting (such as biharmonic-modified) has taken place. For the case of time-dependent partial differential equations (PDEs), after space discretization has taken place what is left is an ordinary differential equation (ODE) system in the time variable which can be numerically integrated by IMEX methods.\\

These methods can also be seen as multi-step methods and they have been
successfully applied to a number of problems including the Navier-Stokes
equations, systems of conservation laws and more recently to
reaction-diffusion equations \cite{splitNVE,shu,Song2016}. 

\par
We now use a general form for 
IMEX methods \cite{ascher1995implicit,Song2016} 
to discretize \eqref{odeF}
as an $s$-stage scheme:
\begin{subequations}
\lbl{imexgeneral}
\begin{eqnarray}
U_{(0)}&=&U_{n}, \\
U_{(i)}&=&\sum_{j=1}^{i} a_{i,j}U_{(j-1)}+h  \sum_{j=1}^{i}
\left( b_{i,j}
\Fe(U_{(j-1)})+c_{i,j} \Fi(U_{(j)})\right)            , \;\;\; i=1,2,\cdots, s, \\
U_{n+1}&=&U_{(s)},
\end{eqnarray}
\end{subequations}
where $s$ indicates the number of stages in the scheme and the coefficients for the scheme are determined based on numerical consistency and order of accuracy 
\cite{shu}.
The $a_{i,j}$ coefficients 
satisfy $\sum_{j=1}^{i}a_{i,j}=1$ to yield consistency in the schemes \cite{shu}.
 The  lower triangular form for $c_{i,j}$ allows for an efficient implementation of the scheme and is used in diagonally implicit Runge-Kutta (DIRK) schemes. 
\par
In this paper, we consider a three-stage ($s=3$) and second-order energy-stable method \cite{Song2016}
\begin{align}
\lbl{imex}
\begin{split}
U_{(0)}&= U_{n}, \hskip3in \mbox{(BHM-IMEX1)}\\
U_{(1)}&= U_{(0)}+h \left( \,  \Fe(U_{(0)})+\Fi(U_{(1)})\,  \right), \\
U_{(2)}&= {\textstyle \frac{3}{2}} U_{(0)} - \half U_{(1)}  + h \left(\,\half
\Fe(U_{(1)})+\half\Fi(U_{(2)}) \, \right), \\
U_{(3)}&=U_{(2)}+ h\left(  \, \Fe(U_{(2)})+\Fi(U_{(3)}) \, \right), \\
U_{n+1}&= U_{(3)}, 
\end{split}
\end{align}
and a second two-stage IMEX scheme (BHM-IMEX2) which is also second-order method and can be expressed in the form \cite{ascher1995implicit} 
\begin{align}
\lbl{imex2}
\begin{split}
U_{(0)} &= U_{n},\hskip3in \mbox{(BHM-IMEX2)}\\
U_{(1)} &=U_{(0)}+h  \left(\, \gamma   \Fe(U_{(0)})+\gamma\Fi(U_{(1)})  \,\right),\\
U_{(2)} &= U_{(0)} + h\left( \,
\delta \Fe(U_{(0)}) +(1-\delta) \Fe(U_{(1)})+ (1-\gamma) \Fi(U_{(1)}) + \gamma \Fi(U_{(2)})\,  \right) ,\\
U_{n+1}&=U_{(2)},
\end{split}
\end{align}
where $\gamma=1-1/\sqrt{2}$ and $\delta= -1/\sqrt{2}$.
The above IMEX time-stepping formulations are motivated from their
successful implementation in Cahn-Hilliard problems with constant mobility \cite{Song2016} and concentration-dependent mobility \cite{HectorC2013}. Both formulations proposed in \cite{Song2016,HectorC2013} are different from what is presented in this paper since the basis for their methods relied on either a convexity splitting approach \cite{Song2016} or a variable mobility product rule expansion \cite{HectorC2013}.

\section{Numerical Tests}
\par
We now examine the performance of these numerical methods for solving Cahn-Hilliard and thin film problems.
The test examples that we present in the spirit of benchmark problems \cite{wettonwise2019} will be used to illustrate the accuracy and stability of the numerical schemes and the dependence on the splitting parameters. Further examples are also presented in the Supplementary Material.

\subsection{Tests of accuracy}
\lbl{sec:accuracy}

First, we present convergence plots for $h\to 0$ to show the order of accuracy of the schemes. Particular attention will be focused on how the fourth-order splitting parameter $\Mb$ influences the error. 
The $\Mh$ term in $\Fi$ also contributes to this, but will be discussed in more detail in section~\ref{sec:energy} in connection with energy stability.

\par
The values used for $\Mb$ can be expected to depend on the range spanned by the mobility function; it will be convenient to define the maximum value of the mobility exhibited by the solution at a given time,
\begin{equation}
\MB(t_n)=\max_{{\bf x}\in \Omega} |\cM(U_n)|\, .
\lbl{MB}
\end{equation}
Indeed, some papers \cite{bertozzi,Glasner2003,barrettblowey,hongkai2003} have suggested that the splitting parameter should satisfy 
$\Mb \ge \MB_{\max}$ where
\begin{equation}
\MB_{\max}=\max_{0\le t_n\le T} \MB(t_n),
\lbl{MBmax}
\end{equation}
namely, $\Mb$ should exceed the mobility everywhere over the whole simulation.
However, we will show that this lower bound is not sharp, and improving this estimate can be important for improving the accuracy of simulations. Indeed, some papers have pointed to the accuracy of methods using splitting as being a significant concern, up to the level of qualitatively changing solution trends for large-time behaviors \cite{christlieb2014high}.

\subsubsection{Test problem 1: accuracy with fixed splitting parameters}   
\par
We consider a typical problem for the dewetting thin film equation (\ref{tfe}), 
\lbl{testprob1}
\begin{equation}
\lbl{tf}
{\partial u\over \partial t} = \nabla\cdot
\left[u^3 \nabla \left(
{\textstyle \frac{\veps^2}{u^3}-\frac{\veps^3}{u^4}}-
 \Delta u\right)\right]\qquad \mbox{on $\Omega=[0,12\pi]^2$}\,,
\end{equation}
with $\veps=0.1$ and periodic boundary conditions. In terms of forms \eqref{CHVMeqn} and \eqref{nlfo2}, this corresponds to 
$\cM(u)=u^3$, $\cG(u)= - 4\veps^3/u - 3\veps^2 \ln(u)$, and
$\cW(u)= -\veps^2 u^{-2}/2 + \veps^3u^{-3}/3$.
For initial conditions, we use a nearly-uniform film with a small-amplitude
spatial oscillation,
$u_0(x,y)=0.35 + 0.1\cos(x+y)$ at $t=0$.
The initial stage of the dynamics, analogous to spinodal decomposition in
the Cahn-Hilliard equation, can be estimated from the prediction of 
linearized dynamics for small perturbation about the unstable constant
state $\bar{u}=0.35$. To ensure that our simulations capture the influences of
the nonlinearities of \eqref{tf}, we evolve the solution to time $t=100$ to yield a large-amplitude nontrivial
initial condition. Testing of numerical convergence then consisted of
evolving from this state, $U_0(x,y)$, 
for one unit of time to the final state $U_*(x,y,T)$ at $T=1$, with different equally-spaced timesteps.
The spatial discretization is implemented  using Fourier pseudo-spectral methods
with $256 \times 256$ modes.


\begin{figure}
\begin{center}
\mbox{
\includegraphics[scale=0.53]{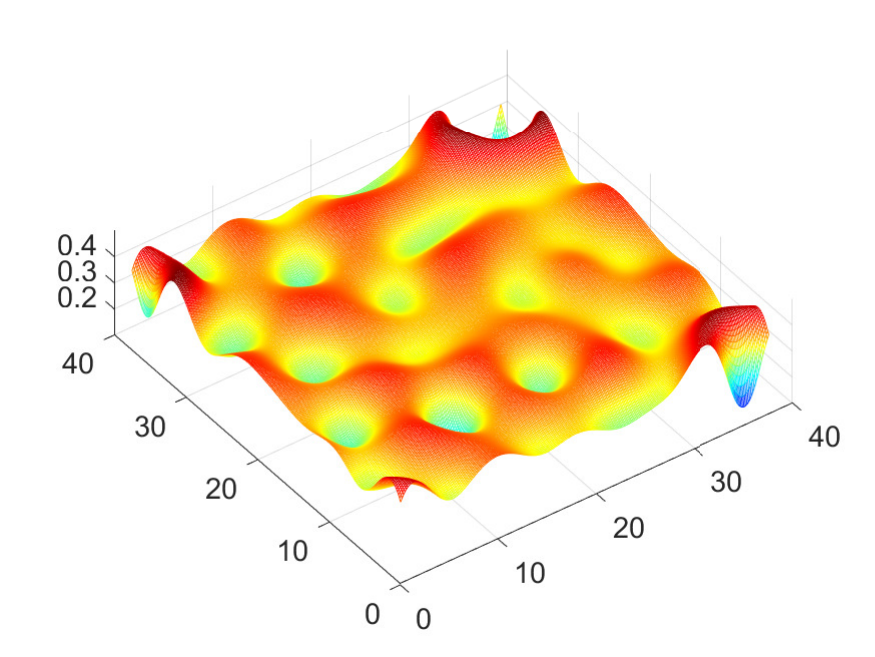}\quad \includegraphics[scale=0.53]{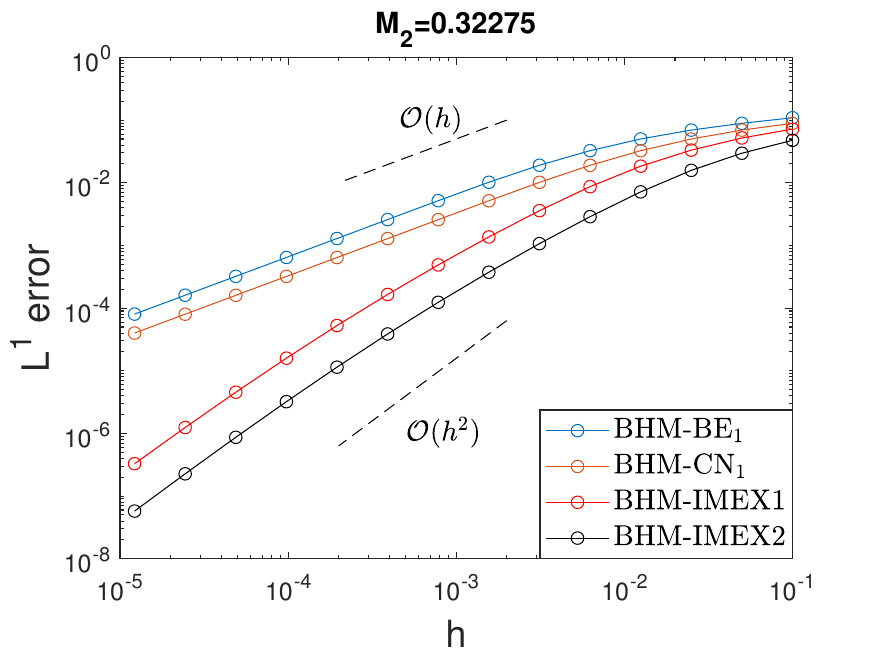}
}
\end{center}
\caption{Thin film test problem 1: (left) The final state of the reference solution $U_*({\bf x},T)$ at $T=1$. (right) Convergence plots for  the four non-iterative schemes
(\ref{BHMBE}, \ref{BHMCN}, \ref{imex}, \ref{imex2}), with splitting parameters $\Mb=0.32275, \Mh=0$.
}
\lbl{fig1}
\end{figure}

\par
An accurate numerical reference solution $U_*({\bf x},t)$ was constructed using
Richardson extrapolation with a very small time-step, $h=10^{-7}$, 
and the errors were computed in terms of the  $L^1$ norm by comparing
solutions for various time increments against the reference solution at
final time $T=1,$ 
\begin{equation}
\mbox{Error}(h) = \int_{\Omega} |U({\bf x},T;h) - U_*({\bf x},T)|\, d{\bf x}.
\end{equation}
The solution was observed to cover the range of mobilities, $0.1227 < \MB(t)< 0.1242$ on $0\le t\le T$.
\par
Figure~\ref{fig1}(right) shows a comparison of the four basic
(non-iterative) time-stepping methods we
consider: \eqref{imex}, \eqref{imex2}, along with \eqref{BHMBE} and 
\eqref{BHMCN} both with $J=1$.
The two single-step methods exhibit first-order accuracy with \eqref{BHMCN}
having a smaller scaling coefficient. The two IMEX methods yield
second-order accuracy with IMEX2 producing consistently smaller errors.
We note that all of these methods were tested with fixed values for the splitting constants of
$\Mb=0.32275$ and $\Mh=0$; this value of $\Mb$ satisfies
$\Mb> \MB_{\max}$.
We will return to issues connected
with the choice of these parameters below.
\par
Figure~\ref{fig2}(left) compares the accuracy of IMEX1, IMEX2 and BHM-BE$_1$ with the iterative extension of BHM,  BHM-BE$_J$ \eqref{BHMBE}, for 
$J=2, 4, 8$. Increasing the number of iterations decreases the error, but increases the
computational work per timestep. Using pseudo-spectral methods, the linear implicit terms are
straightforward, so the main computational workload is attributable to each evaluation of the
nonlinear explicit operator $\Fe(U)$. The BHM-BE$_J$ method has $J$ evaluations of $\Fe$ per timestep compared
with 3 evaluations per timestep for both IMEX1 and IMEX2.
After two iterations, BHM-BE$_2$ has an error
that is only slightly larger than IMEX1 with somewhat lower computational work
involved. Increasing the number of iterations to $J=4$ to compete with IMEX2 shifts the speed advantage to IMEX2.
An unusual aspect of the error curves for the BHM-BE$_J$ iterative
methods is that they appear to scale like $Ch^2$, like the curves for the second-order IMEX methods.
However, this is not actually the asymptotic behavior for the BHM-BE$_J$ methods for 
$h\to 0$; we will see that that it only holds for an intermediate
range of time-step sizes. Increasing the number of iterations lowers the error at each fixed value of $h$. The ``bend'' in the error curve for BHM-BE$_8$ on $10^{-5}\le h\le 10^{-4}$ recovers the first-order accuracy expected for Backward Euler methods.
We will discuss how this behavior is influenced by 
the implicit/explicit splitting parameter $\Mb$ next.

\begin{figure}
\begin{center}
\mbox{
\includegraphics[scale=0.53]{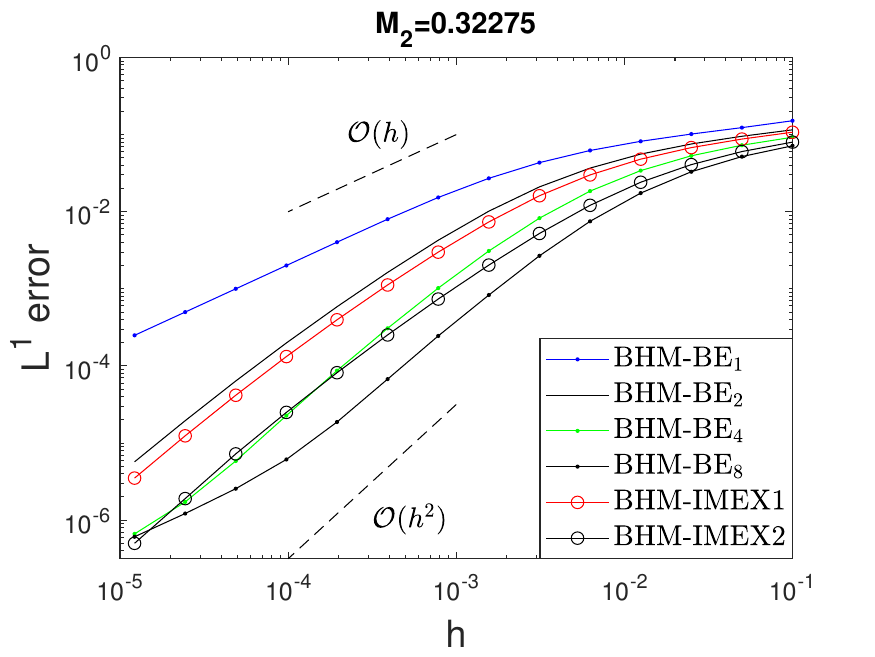}\quad 
\includegraphics[scale=0.53]{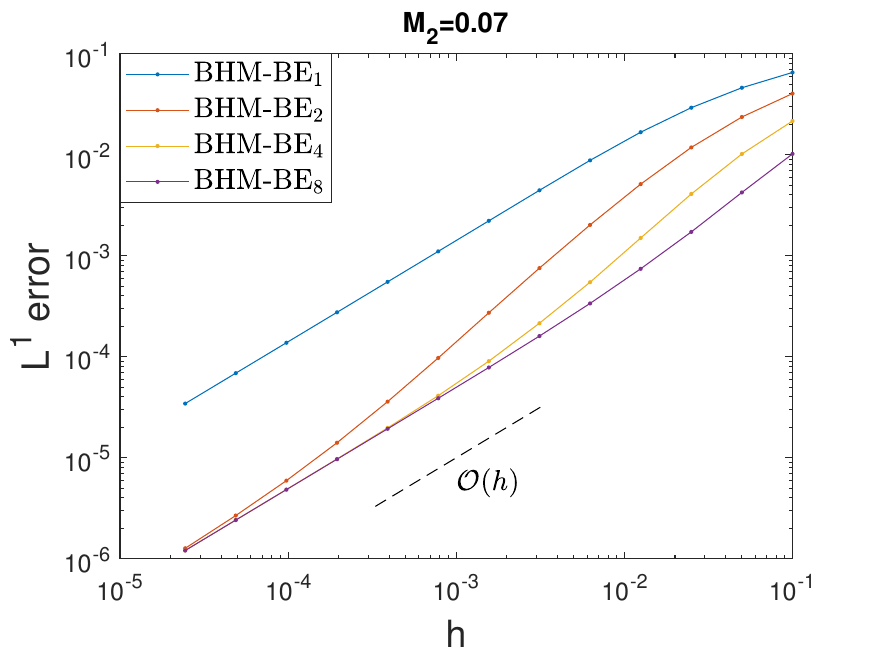}
}
\end{center}
\caption{ 
Test problem 1 (continued): 
(left) Comparison of the IMEX schemes with the
iterative BHM-BE$_J$ scheme for $J=1,2,4,8$, all with $\Mb=0.32275$ and $\Mh=0$. 
(right) Error curves for the BHM-BE$_J$ iterative methods with
splitting parameters $\Mb=0.07, \Mh=0$. 
}
\lbl{fig2}
\end{figure}

\par
In Figure~\ref{fig2}(right)
the first-order accuracy of the BHM-BE$_J$ iterative methods is made more clear through the use of a smaller value of the splitting parameter, here $\Mb=0.07$
compared to $\Mb=0.32275$
used in the earlier figure. Here increasing the number of iterations exposes the first-order rate of convergence over a wider range of timesteps. We conjecture that for large numbers of iterations ($J\ge 4$) the error curve could approach the error produced by the fully implicit (nonlinear) backward-Euler method. 
The reduction in the error relative to the non-iterative BE$_1$ method is seen to be over a factor of ten.
The numerical results suggest that more careful analysis beyond \eqref{TruncErr} would show the truncation error to scale like $||\tau||\sim C_1(\Mb;J)h + C_2(\Mb;J)h^2$ with a relatively large $C_2$ coefficient.

\par
In simulations of each of the numerical schemes,
 we observed that at any fixed $h$, decreasing $\Mb$ decreases the error. 
There must be a positive lower bound on $\Mb$ for stable simulations.   
If $\Mb$ were set to zero, then the methods reduce
to explicit methods that have severe limitations on the maximum timestep, $h=O(\Delta x^4)$ for conditional stability for stiff fourth-order parabolic PDEs.
Figure~\ref{fig:M1effect} shows error curves for three methods (BHM-BE$_1$, IMEX1, IMEX2) over a range of values for $\Mb$ with a fixed value of the timestep, $h=0.125$. 
Above a critical minimum value of the splitting parameter, $\Mb> \Mb^*$, the errors are monotone decreasing for each scheme.  
Below $\Mb^*$ numerical instabilities producing large errors were observed.
Estimates for these critical values were found as 
$M^*_{2,\mathrm{BE1}}\approx 0.0528$, $M^*_{2,\mathrm{IMEX1}}\approx 0.0628$, and $M^*_{2,\mathrm{IMEX2}}\approx 0.120$.
These values are specific to this test problem but to interpret them in a broader context, we can re-cast them relative to \eqref{MBmax} using
\begin{equation}
\Mb=  \overline{\alpha} \;\MB_{\max}
\lbl{aMmax}
\end{equation}
to give  
$\overline{\alpha}^{\,*}_{\mathrm{BE1}}\approx 0.43$, $\overline{\alpha}^{\,*}_{\mathrm{IMEX1}}\approx 0.51$, and $\overline{\alpha}^{\,*}_{\mathrm{IMEX2}}\approx 0.97$.
These estimates are based on a single fixed timestep size, in section~\ref{sec:energy} we will re-examine this using a different approach.
Noting that all of these are below the case of $\overline{\alpha}=1$ given by \eqref{MBmax}, we observe that it is possible to lower $\Mb$ to thus gain more accurate solutions.
Some previous papers have offered results for single-step schemes suggesting $\overline{\alpha}\ge 1/2$ in \eqref{aMmax} is stable, 
for example \cite{smereka,shenyang2010,vollmayr2003,shengdu2010,zhutikare1999,eggers2014}. The value obtained for $\overline{\alpha}^{\,*}_{\mathrm{BE1}}$ in Fig.~\ref{fig:M1effect} being slightly below the stability criterion is unexpected but may be due to the results being carried out at a moderate-sized timestep.
Other papers have suggested that using $\overline{\alpha} > 1$ could yield possible instabilities \cite{he2007large}. Returning to the convergence results shown in Figure~\ref{fig2}, in terms of \eqref{aMmax} the BHM-BE simulations shown there were done with $\overline{\alpha}_{\mathrm{BE}}\approx 2.6$ and
$\overline{\alpha}_{\mathrm{BE}}\approx 0.57$ respectively. 

\par 
Further simulations on this test problem showed that increasing the second-order splitting parameter from 
$\Mh=0$ generally increased the error.  This will be revisited in section~\ref{testprob4} where it is shown that $\Mh$ can improve results from the schemes but must be selected carefully.

\begin{figure}
\begin{center}
\mbox{
\includegraphics[scale=0.57]{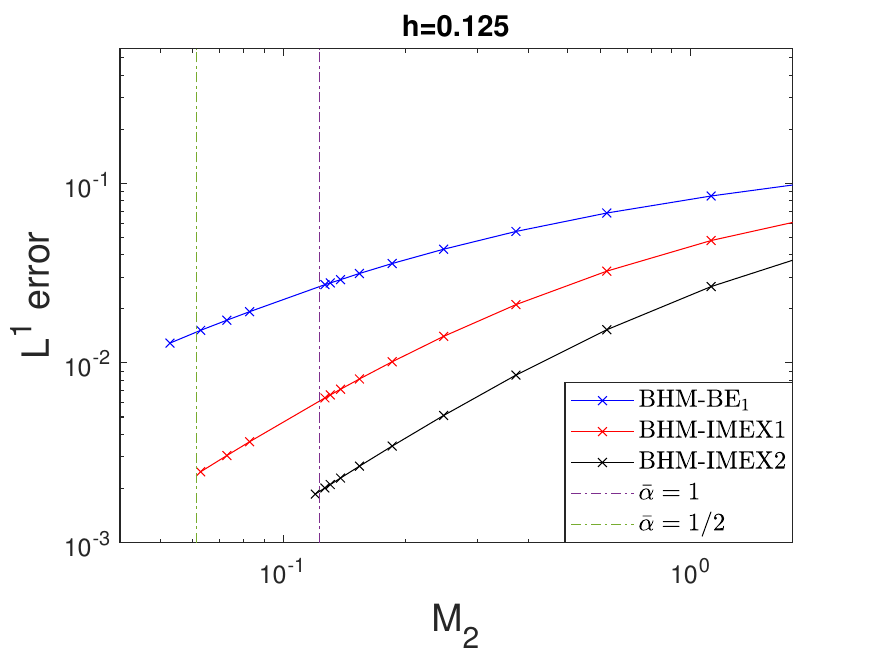}
}
\end{center}
\caption{
Test problem 1 (continued): 
The error for the time-stepping methods at fixed $h=0.125$ over a
range of values for the $\Mb$ splitting parameter. The curves terminate at minimum values $\Mb^*$ below which simulations went unstable. Vertical lines at fixed values for the $\overline{\alpha}$ ratio in \eqref{aMmax} are shown for reference. }
\lbl{fig:M1effect}
\end{figure}


\subsubsection{Test problem 2: static vs.\ dynamic splitting} 
We briefly consider a second test problem to more carefully examine how choices for splitting parameters influence the accuracy of these methods.
\par
A common approach for constructing test problems is by adding an inhomogeneous forcing term to the PDE to make a given function an exact solution, see for example \cite{yanwise2018} for a test problem for the forced Cahn-Hilliard equation. 
For the thin film equation \eqref{tf}, the function
$u_{\mathrm{exact}}(x,y,t)=0.3+0.1\sin(x)\sin(y)e^{t/2}$
can made to be an exact solution by inserting it into equation \eqref{tf} and obtaining the corresponding
forcing term as $\widetilde f(x,y,t)\equiv \partial {u_{\mathrm{exact}}}/\partial t-F(u_{\mathrm{exact}})$. This gives the following modified or forced problem,
\begin{equation}
\lbl{fCH}
{\partial u\over \partial t} =\nabla\cdot
\left[u^3 \nabla \left(
{\textstyle \frac{\veps^2}{u^3}(1-\frac{\veps}{u})}-
 \Delta u\right)\right]+\widetilde f(x,y,t)\,,
\end{equation}
beginning with $u=u_{\mathrm{exact}}(x,y,0)$ as initial data.
\par
Figure~\ref{fig:CNforced}(left) shows error curves for the iterative BHM-CN$_J$ \eqref{BHMCN} methods with $J=1,2,4,8$ and $\Mb=0.125$, $\Mh=0$. Without iteration, the error for the BHM-CN$_1$ schemes shows only first-order accuracy. However, adding even one iteration makes the method approach second-order accuracy ($J\ge 2$), and like the backward Euler results, we expect that with sufficient iterations per time-step, this method may converge to the accuracy of the corresponding nonlinear implicit Crank-Nicolson scheme. In terms of computational workload, the CN$_J$ method requires $J+1$ evaluations of $\Fe$ per timestep vs.\ 3 evaluations per timestep for the two multi-step IMEX methods. In terms of accuracy considerations alone, the CN$_2$ becomes competitive with the IMEX schemes producing slightly better errors than IMEX1 and just under-performing the IMEX2 errors. For higher numbers of iterations CN$_J$, $J=4,8$ generates smaller errors compared to both IMEX schemes but for many iterations the CN$_J$ scheme becomes more computationally expensive.

\begin{figure}
\begin{center}
\mbox{
\includegraphics[scale=0.57]{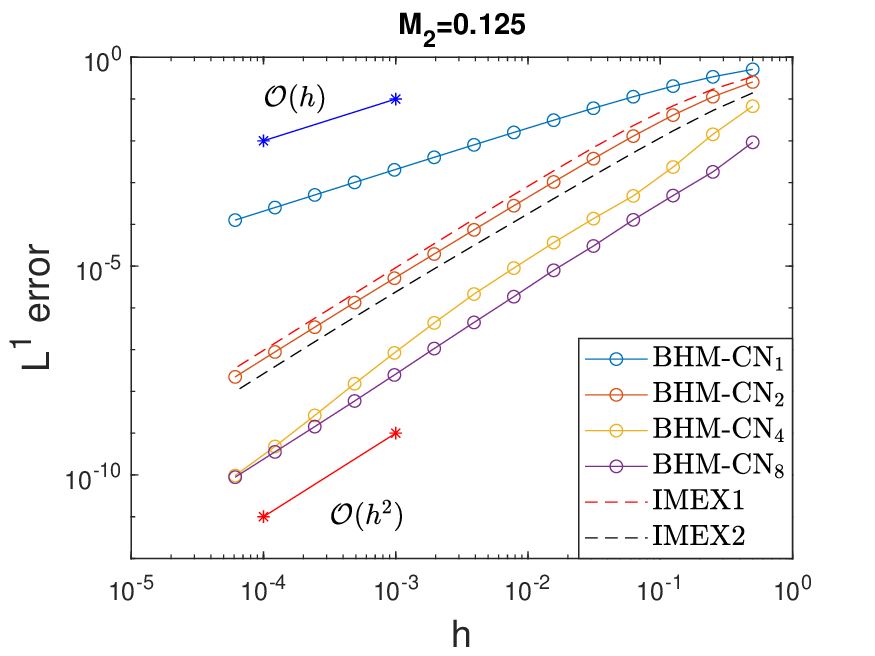}
\includegraphics[scale=0.57]{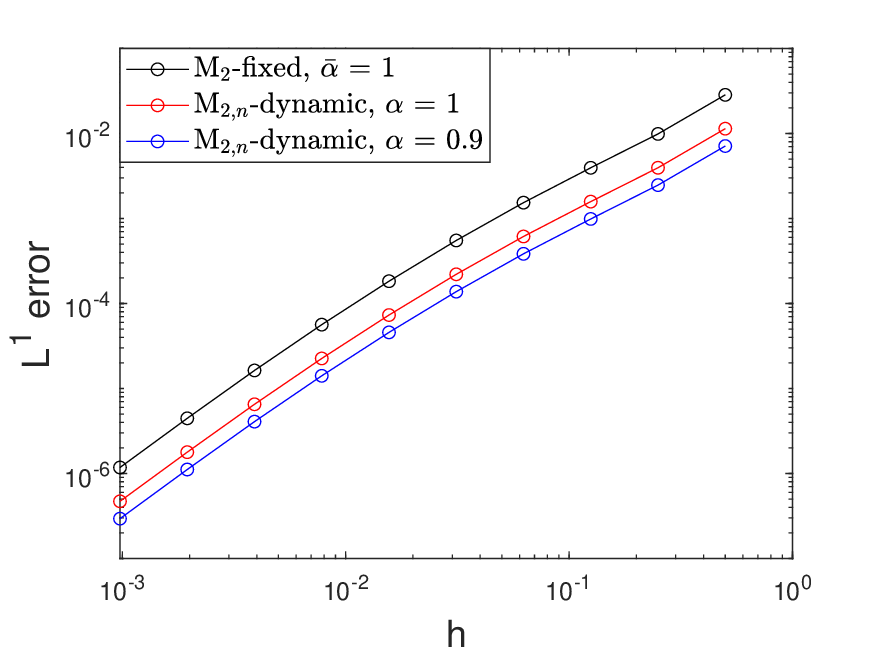}
}
\end{center}
\caption{Test problem 2: (left) Error curves for the  BHM-CN$_J$ iterative methods with $J=1,2,4,8$ and the two IMEX schemes, 
(right) Error curves for IMEX2 using fixed splitting parameter $\Mb=\MB_{\max}$ compared with dynamic splitting with  $M_{2,n}=\MB(t_n)$ and $M_{2,n}=0.9\MB(t_n)$.}
\lbl{fig:CNforced}
\end{figure}

\par
In this test problem, since the solution is growing, the mobility is likewise monotone increasing, having
$0.065< \MB(t)< 0.125$ on $0\le t\le T$ with $T=1.4$. 
Figure~\ref{fig:CNforced}(right) compares the accuracy of the IMEX2 scheme with two different choices for the splitting parameter: (i) a fixed value $\Mb=\MB_{\max}=0.125$ as in Fig.~\ref{fig:CNforced}(left) (i.e. \eqref{aMmax} with $\overline{\alpha}=1$), vs.\ (ii) selecting $\Mb$ to be a dynamic parameter with a different value at each timestep,
\begin{equation}
M_{2,n} =\alpha \;\MB(t_n),
\lbl{M1n}
\end{equation}
used here with $\alpha=1$ and $\alpha=0.9$. 
This kind of fixed-ratio choice for the dynamic fourth-order splitting parameter was also used in \cite{bertozzi}.
Relative to static splitting, dynamic splitting with the same ratio ($\alpha=\overline{\alpha}$) yields smaller values for $\Mb$ all times $t_n$ when the global maximum for $M(U)$ has not been attained, and hence will yield smaller errors at almost all steps.
Figure~\ref{fig:CNforced}(right) shows that using a dynamic $M_{2,n}$ can improve the accuracy by a factor of two and can be further improved for lower values of $\alpha$ down to the lower-bound $\alpha>\alpha^*$ for yield stable computations for the scheme being used.
\par 
In the next section we will consider the stability of computations with different values for the $\Mb, \Mh$ splitting parameters.

\subsection{Tests of energy stability}
\lbl{sec:energy}
\par
Since the Cahn-Hilliard and thin film PDE's both have monotone dissipated energy \eqref{dEdt}, discrete-time energy stability (also called gradient stability \cite{bertozzi})  of the numerical schemes is necessary for consistency of simulation results,
\begin{equation}
\mathcal{E}(U_{n+1})\le \mathcal{E}(U_n)\qquad \forall n\,.
\lbl{Etest}
\end{equation}
Determining if some classes of implicit-explicit numerical schemes are unconditionally energy stable for all timesteps, $h>0$, has received a lot of attention in the numerical analysis literature \cite{rosales,christlieb2013unconditionally}. The interest in large timesteps is motivated by the use of adaptive timestepping that becomes very relevant for practical computations for these PDE models since their long-time behaviors can exhibit very slow dynamics.

\par
Many papers have presented analytical estimates for conditions on energy stability for various numerical methods \cite{christlieb2013unconditionally,Song2016}; here we carry out direct numerical tests to see if
\eqref{Etest} holds for our schemes.
For each test, simulations of the PDE will be run over a range of numerical parameters, the timestep $h$ and one of the splitting parameters ($\Mb$ or $\alpha$ or $\Mh$) and if \eqref{Etest} is satisfied for the entire simulation will be recorded in a plot of the parameter plane. For the figures in this section, success will be marked in yellow, while if \eqref{Etest} fails at any timestep the simulation will be shown in blue (e.g. see Fig.~\ref{fig:CS_Stability}). Our results suggest there are fairly robust dividing curves separating stable from unstable parameter ranges. These tests are not definitive, but their results will be very suggestive of parameter ranges where unconditional stability of the numerical methods could be potentially provable.
Starting from a fixed choice of initial condition, the PDE is simulated using one of the schemes to evolve the solution for a fixed number of time-steps, $N$, independent of the size of $h$. 
{ Further remarks on comparisons with alternative ways to test \eqref{Etest} are given in the Supplementary Material.}
We will present two test problems to show how energy stability can depend
on the $\Mh$ or $\Mb$ parameters separately and where $\Mh, \Mb$ must be considered together.

\begin{figure}         
\centerline{
(b)\includegraphics[scale=0.53]{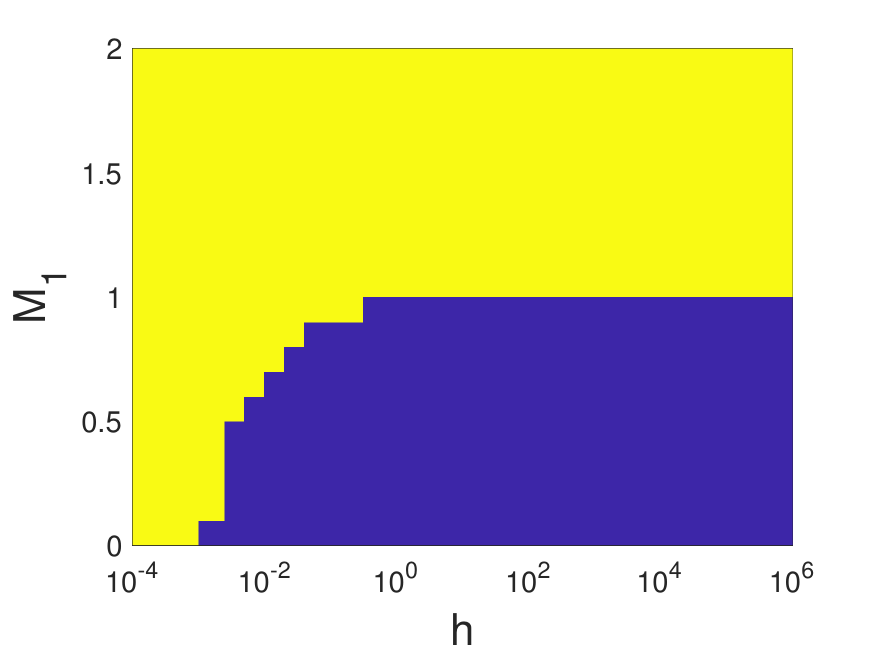}
}
\caption{Test problem 3: Energy stability of numerical solutions of the Cahn-Hilliard equation \eqref{CHclassic} using the BHM-BE$_1$ method. Simulations with $(h, \Mh)$ parameters that are energy-stable are marked yellow, unstable parameters are shown in blue. }
\lbl{fig:CS_Stability}
\end{figure}


\subsubsection{Test problem 3: energy stability for the Cahn-Hilliard equation}
\lbl{cchEnergy}
The Cahn-Hilliard equation \eqref{CHclassic} 
with $\epsilon=0.02$ on $\Omega=[0,\pi]^2$, and a specified initial condition was used as a computational test problem in \cite{glasnerOrizaga}. Matching this PDE with form \eqref{nlfo2} gives $\cM(u)=\epsilon^2$ and $\cG(u)=u^3-u$. 
Since the mobility is a constant, setting $\Mb=\epsilon^2$ gives the full fourth-order operator exactly.
For the second order operator, we write $\cG_-(u)=u^3-u-\Mh u$ and $\cG_+(u)=\Mh u$.  The convexity splitting condition from \eqref{eyre}  yields that $\cG_-'(u)=3u^2-(1+\Mh)<0$ for the range of values for the solution, nominally assumed to be $-1\le u\le 1$ suggests that $\Mh>2$ should be used. The numerical method used in \cite{glasnerOrizaga} is equivalent to the BHM-BE$_1$ scheme \eqref{schemeBH} with $\Mb=\epsilon^2$ and for specified values of $\Mh$. Interestingly, in \cite{glasnerOrizaga}, the value $\Mh=3/2$ (below the convexity splitting threshold) was used and produced accurate, energy-stable results. 
\par
Figure~\ref{fig:CS_Stability}(a) shows the energy stability results for this problem for a range of $\Mh$ 
with $N=500$ time-steps for $10^{-4}\le h\le 10^6$. 
These results suggest that the method could be unconditionally energy stable for all $\Mh>1$ \cite{he2007large}.
This supports the use of $\Mh=3/2$ from \cite{glasnerOrizaga} as a stable splitting that would be more accurate than using $\Mh=2$.

\subsubsection{Test problem 4: energy stability for the thin film equation }
\lbl{testprob4}
We now consider the thin film equation \eqref{tf} with $\veps=0.1$ on $\Omega=[0,6\pi]^2$ and $0\le t\le 200$ with initial conditions describing a weakly perturbed film with mean thickness $\bar{u}\approx 0.6$. 
The dynamics of the solution are illustrated in Figure~\ref{fig:IC_Stability2}. In the early stages, the amplitude of perturbations grows and one local minimum dominates to approach near-rupture. The rupture point then grows to form a well-defined ``hole'' with $u\approx \veps$.
The maximum of the mobility is found
to be $\MB_{\max}\approx 0.621$ achieved near $t\approx 200$ having started from $\MB_{\max}\approx 0.216$ at $t=0$ (see Figure~\ref{fig:IC_Stability2}(c)).
We will now see that there are interactions between the second- and fourth-order splitting parameters involved in the energy-stable of the thin film equation.
\par
We will describe the results obtained for the IMEX-1 scheme applied to this problem. We found that all of the methods considered in this paper yield qualitatively similar parameter dependence for their energy stability. The plots in Figure~\ref{fig:TF_Test5} show that for the thin film equation it was not possible to obtain consistent solutions for arbitrarily large time-steps. This is suggestive of the loss of accuracy becoming the limiting factor in simulations that might still be formally unconditionally stable \cite{eggers2014}. Fig.~\ref{fig:TF_Test5}(left) shows that IMEX-1 simulations with $\Mh=0$ and $\Mb$ set by \eqref{M1n} are energy stable for moderate time-steps for $\alpha\ge \alpha^*_{\mathrm{IMEX1}}\approx 0.5$. For $\alpha< \alpha^*$ the limitation on the maximum time-step is comparable to the linear stability bound for explicit methods, $h=O(\Delta x^4)$. Our further energy stability tests suggest the critical fourth-order parameters for the other numerical schemes are $\alpha_{\mathrm{BE1}}\approx 0.5$, $\alpha_{\mathrm{CN1}}\approx 1.0$, and $\alpha_{\mathrm{IMEX2}}\approx 0.85$.
\par
The Cahn-Hilliard and thin film equations can both be expressed in terms of \eqref{nlfo2}, but differences with respect to their second order operators may be very significant.
The second order nonlinearity
$\cG(u)$ for \eqref{tf} can be split similarly to what was done in section~\ref{cchEnergy}. Requiring
convexity of $\cG_-(u)=\cG(u)-\Mh u$ for $u\ge \epsilon$ suggests that the
splitting parameter can be any value $\Mh\ge \epsilon$. However, it appears
that this condition is not necessary for energy stability, as Figure~\ref{fig:TF_Test5}(left) used 
$\Mh=0$ and still showed ranges of energy stability.
\par
Figure~\ref{fig:TF_Test5}(right) shows results for the IMEX1 with $\alpha=1$ for a large range of values for the second-order splitting parameter $\Mh$. Above a certain critical value for $\Mh$, here $\Mh^*\approx 1000$, the influence of the second-order splitting allows for the use of significantly larger timesteps. Above $\Mh^*$ the size of the maximum timestep depends on $\Mh$ and this requires further study. However, for a fixed value of $\Mh\ge\Mh^*$ we have observed that $h_{\max}$ is larger than in Fig.~\ref{fig:TF_Test5}(left) and is generally independent of the value of $\alpha>\alpha^*$ used for the fourth-order splitting. This suggests that the second-order splitting is primarily responsible for the energy stability. This may fit with the results of Duchemin and Eggers \cite{eggers2014} showing that numerical schemes for PDE can be stabilized using lower-order operators.

\begin{figure}
\begin{center}
\mbox{
(a)\includegraphics[scale=0.35]{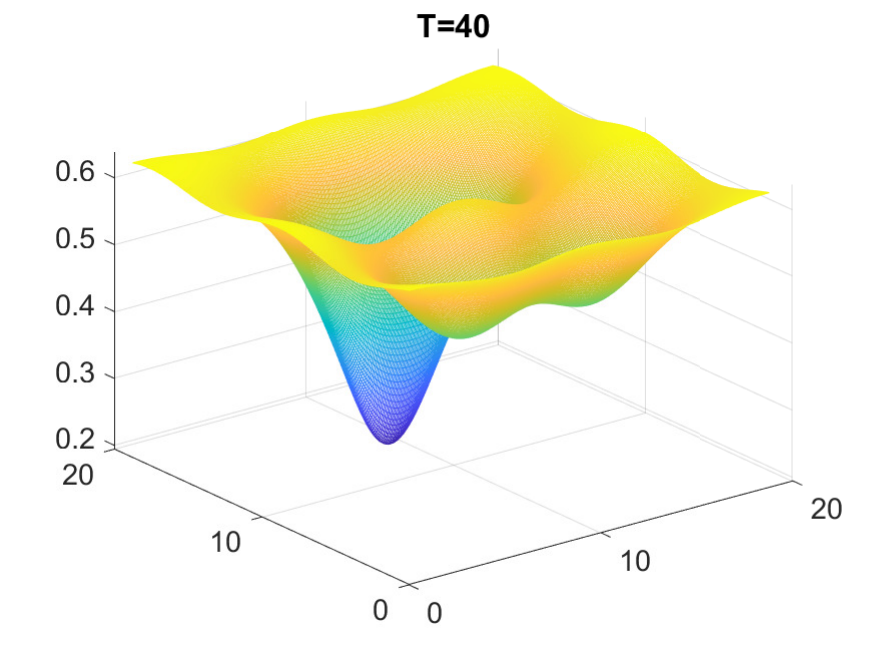}
(b)\includegraphics[scale=0.35]{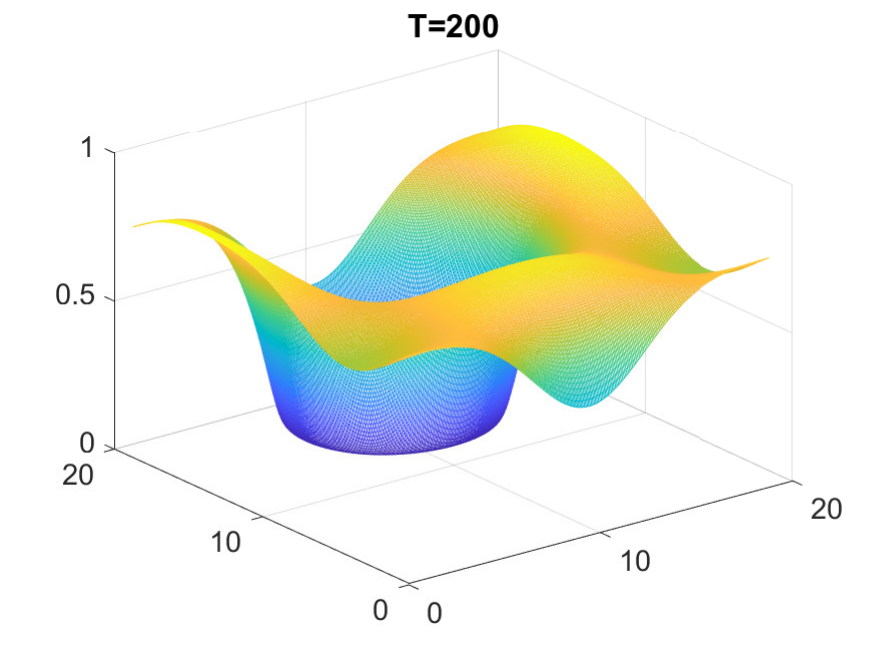}
(c)\includegraphics[scale=0.35]{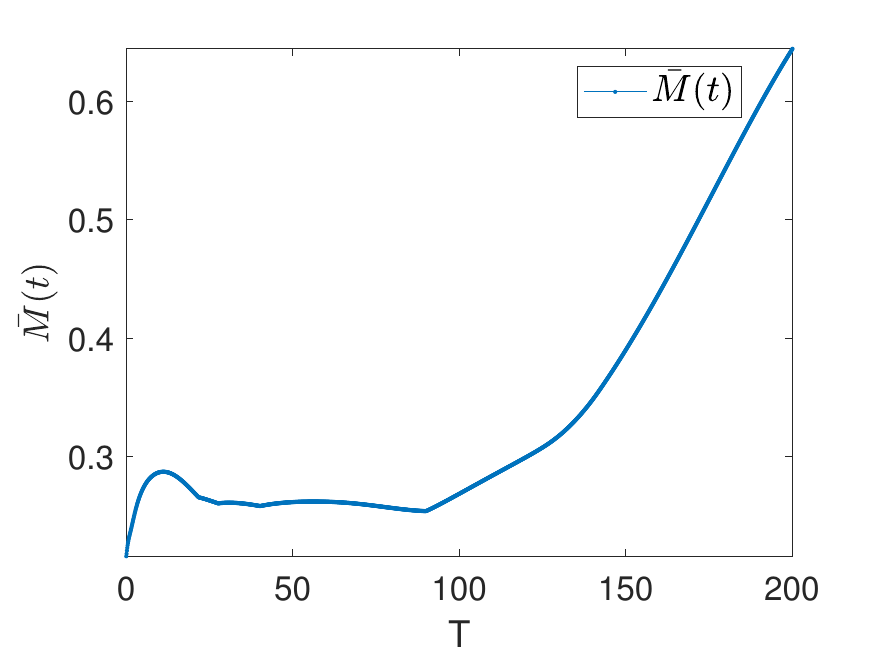}
}
\end{center}
\caption{Test problem 4: Numerical solution of the thin film equation \eqref{tf} using BHM-IMEX2 with $\Mb=1$, $\Mh=0$ and $h=0.001$. 
(a) Solution at time $t=40$. (b) Solution at time $t=200$. (c) Evolution of the maximum of the mobility coefficient, $\MB(t)$, \eqref{MB}.}
%
\lbl{fig:IC_Stability2}
\end{figure}

\begin{figure}
\begin{center}
\mbox{
\includegraphics[scale=0.5]{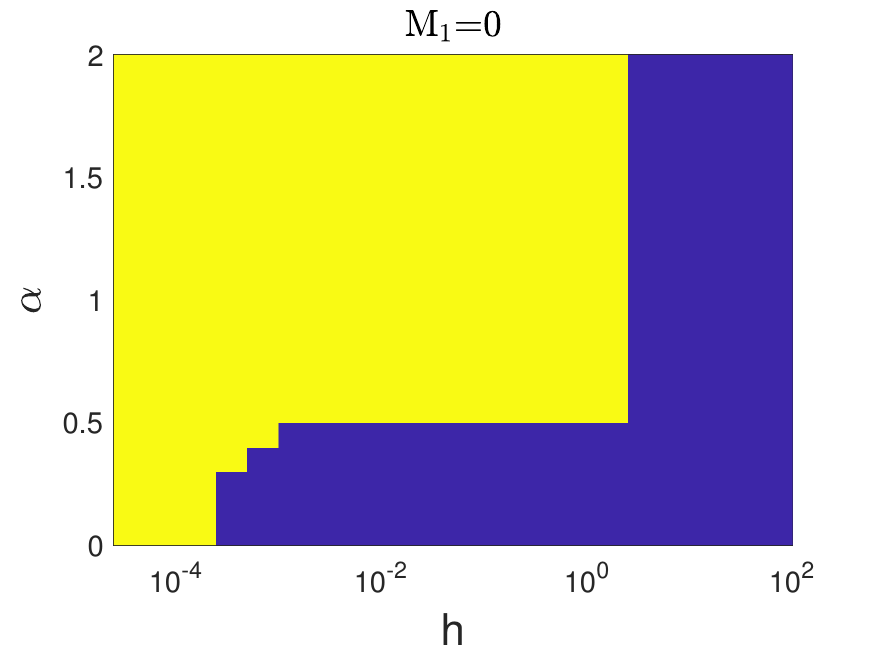}
\includegraphics[scale=0.5]{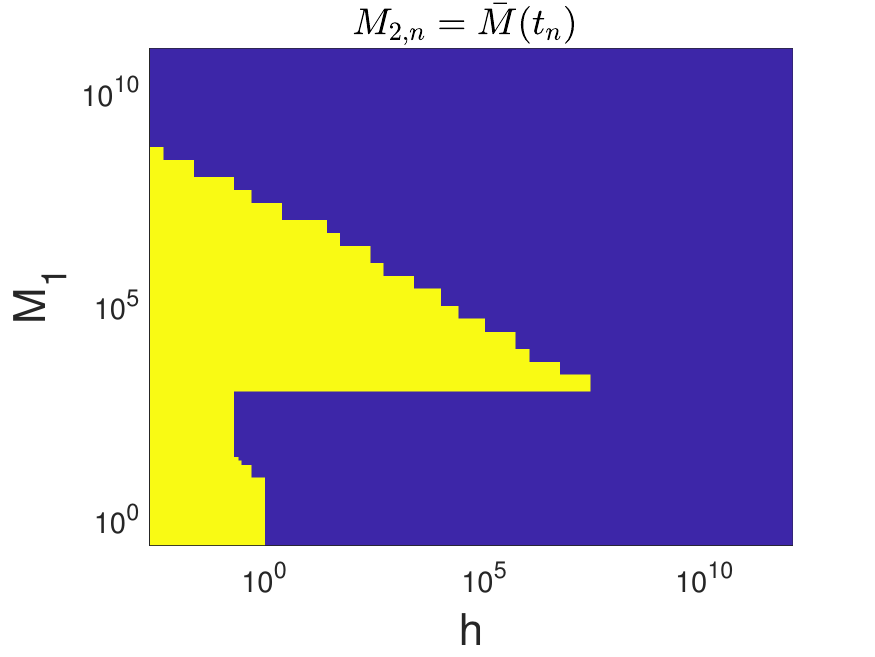}
}
\end{center}
\caption{
Test problem 4: Energy stability plots for the IMEX1 method applied to the thin film equation \eqref{tf}. (left) Dependence of the energy stability on the $\alpha$-ratio for the fourth-order dynamic splitting parameter \eqref{M1n}. (right) Dependence on the second-order splitting parameter $\Mh$ with $\alpha=1$ in \eqref{M1n}.}
\lbl{fig:TF_Test5}
\end{figure}

\par
Our simulations suggest that the methods are not able to obtain true unconditional energy stability (i.e. \eqref{Etest} for unbounded timesteps). This is attributed to the fact that the variable mobility $\cM(u)$ makes the component in the TF equation that needs to be computed implicitly to become non-linear, $ \div\left(\cM(u) \grad \left[\frac{\delta \mathcal{E_+}}{\delta u}\right] \right)$, and this implicit computation is mitigated by partially performing the task with the fourth-order splitting $\Mb$. In order words, the unconditional energy-stability property is not guaranteed since the splitting applied to TF or CH equations with variable mobility becomes a modified version of the original convexity-splitting scheme proposed by Eyre  \cite{eyre1998unconditionally}.
The present energy stability results provide a basic description of the parameter space for the fourth-order and second-order splittings. These results could provide a guideline for choosing energy-stable parameters in a given simulation and help to understand and quantify the advantages and limitations of the currently available splittings. Future investigations should explore and develop new splittings or methods that improve the energy-stability property with respect to the choices in the parameter space. 






\section{Illustrations of computed dynamics}

Dynamics in thin-film problems may take a longer time to exhibit full drop formation, so our schemes must be able to perform accurately and efficiently for longer runs.  We illustrate the performance of  BHM-IMEX2 by solving the tf equation \eqref{tf} with the same initial condition as in section~\ref{testprob1}  and all parameters the same except that we consider a larger box $\Omega =[0, 24\pi]^2$ and $\tf=1250$.
The motivation for choosing a larger domain relies simply in the desire to capture more well-defined drops at the end of the simulation. Using $\Mb=5.0$ and $h=0.1$ the simulation snapshots for  $t=250,500,750$ and $1250$ (left to right) are shown in  Figure~\ref{fig:thdirk}. The method exhibits the expected dynamics of a thin-film undergoing instabilities that result in drop formation. The simulation results at $t=1250$  illustrate well-defined drops that slowly evolve due to mass exchange with remaining drops which is consistent with the physical evolution for the thin-film. For this numerical experiment mass is conserved and the solution preserves the energy-decreasing property.
Similar numerical simulations capturing the correct dynamics were observed for the BHM-BE$_J$, BHM-CN$_J$ and BHM-IMEX1 methods.

We note that while the methods considered in this paper are energy-stable and depending on the accuracy requirements at hand, a timestep smaller than $h=1$ would be recommended to obtain more accurate solutions. However for computational tasks in which accuracy is modest, a large time step such as $h=1$ can serve the purpose of a quicker computation that allows to explore and investigate the overall dewetting dynamics associated with the solution to Eq (\ref{tfe}).

\begin{figure}
\begin{center}
\mbox{
\includegraphics[scale=0.23]{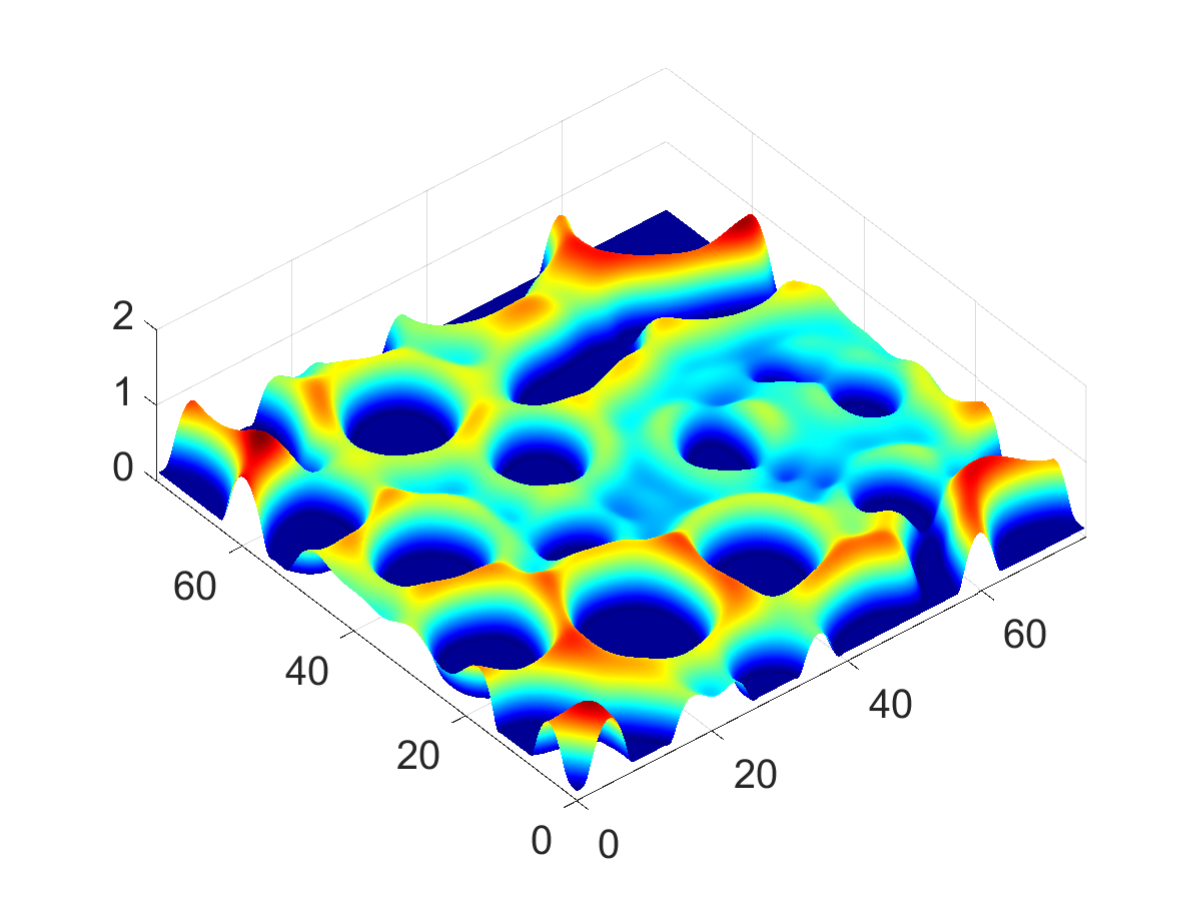}\hskip-0.15in
\includegraphics[scale=0.23]{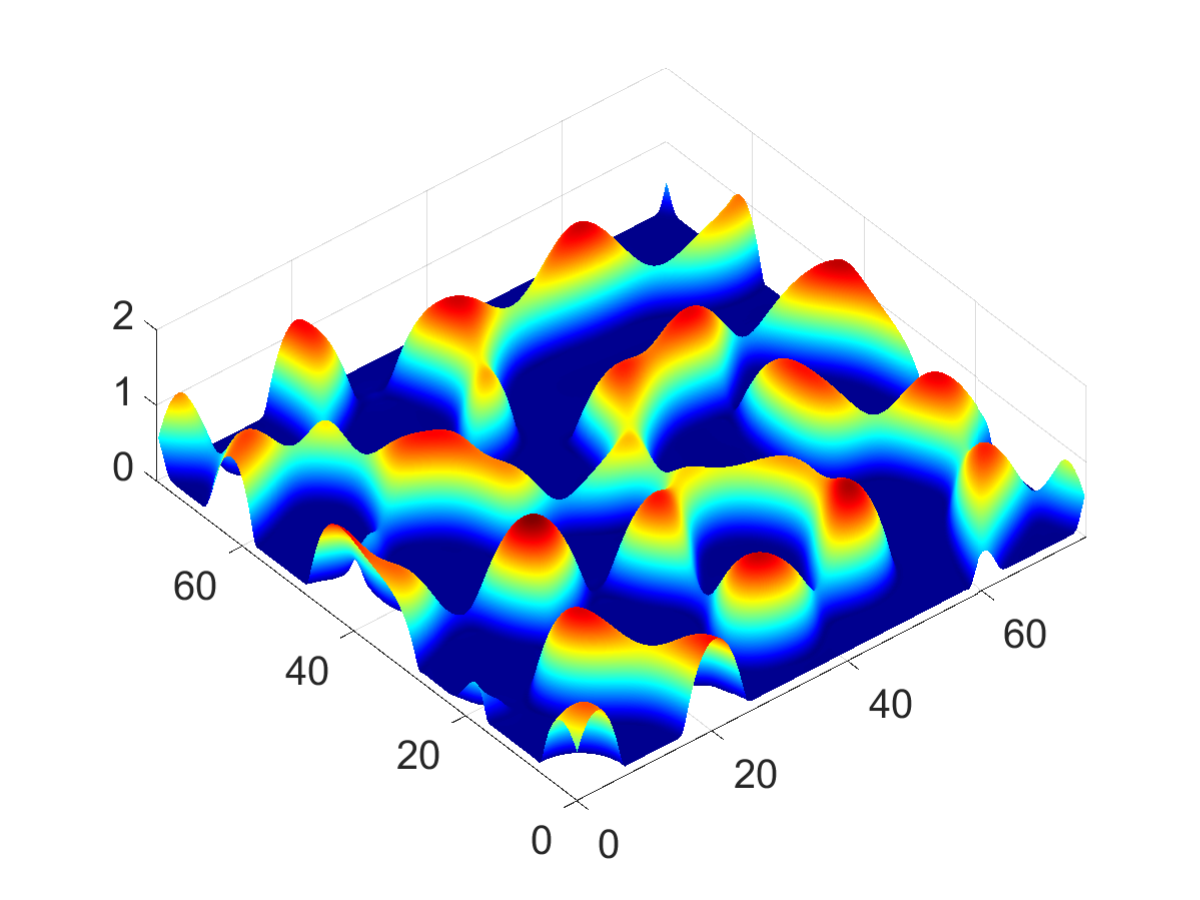}\hskip-0.15in
\includegraphics[scale=0.23]{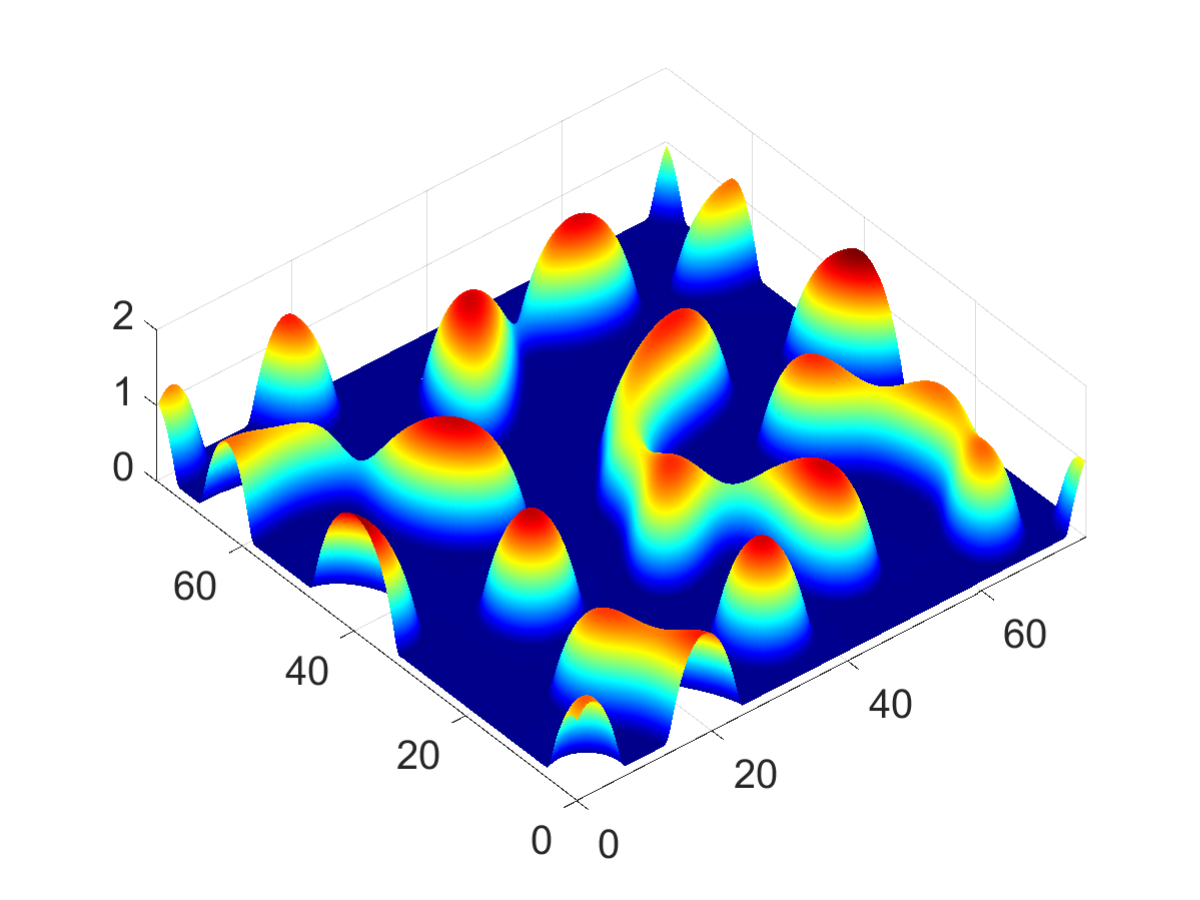}\hskip-0.15in
\includegraphics[scale=0.23]{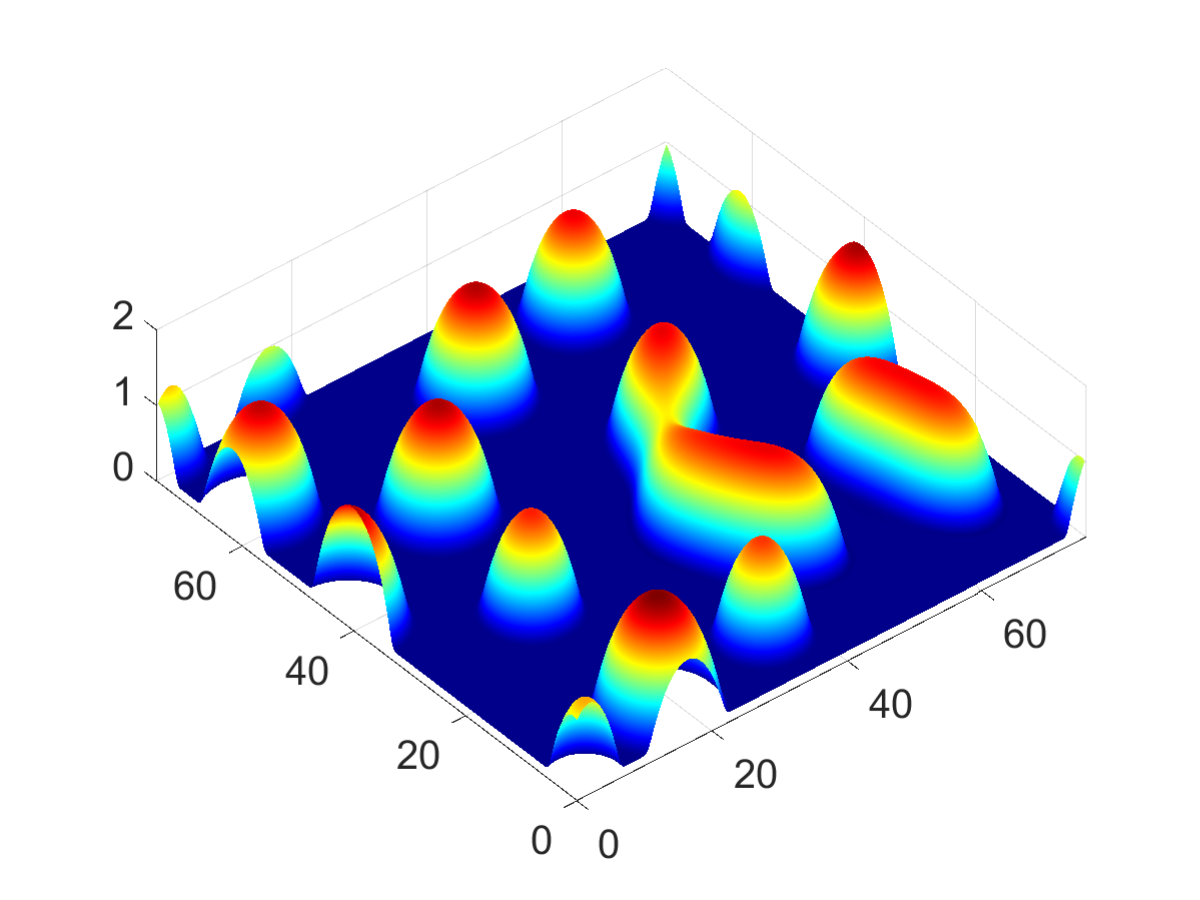}}
\end{center}
\caption{
Evolution of the numerical solution to the thin film equation \eqref{tf} using the BHM-IMEX2 method with  $\Omega=[0,24\pi]^2$, $\Mb=5$, $h=0.1$ with profiles at times $t=250,500,750, 1250$.
}
\lbl{fig:thdirk}
\end{figure}



\begin{figure}
\begin{center}
\mbox{
(a)\includegraphics[scale=0.43]{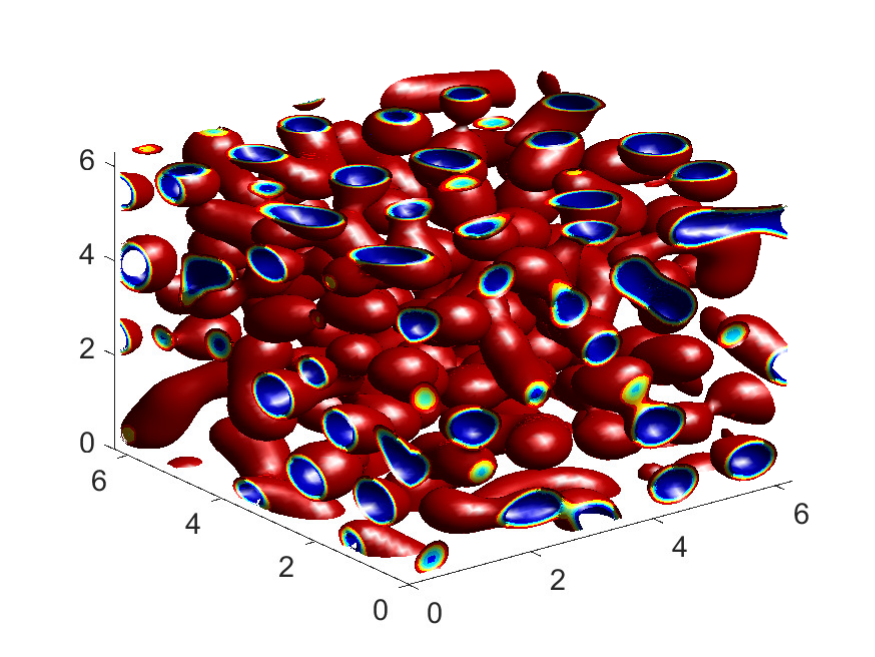}
(b)\includegraphics[scale=0.43]{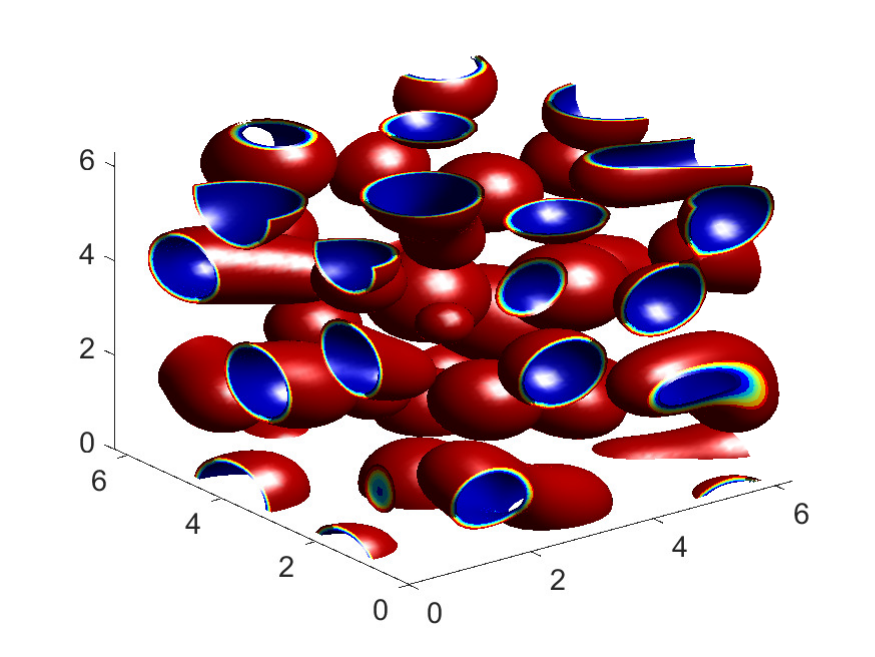}
}
(c)\includegraphics[scale=0.43]{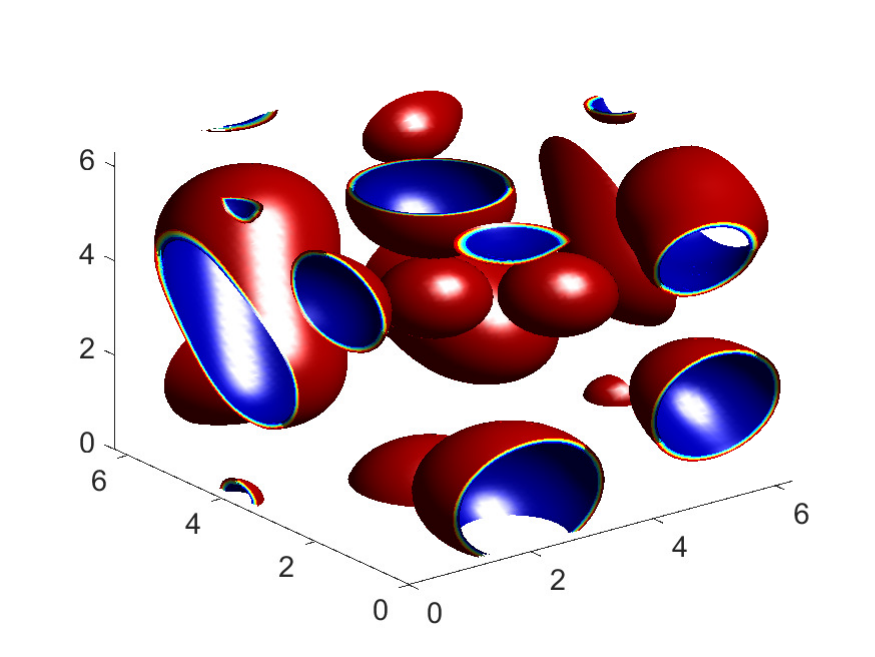}
(d)\includegraphics[scale=0.43]{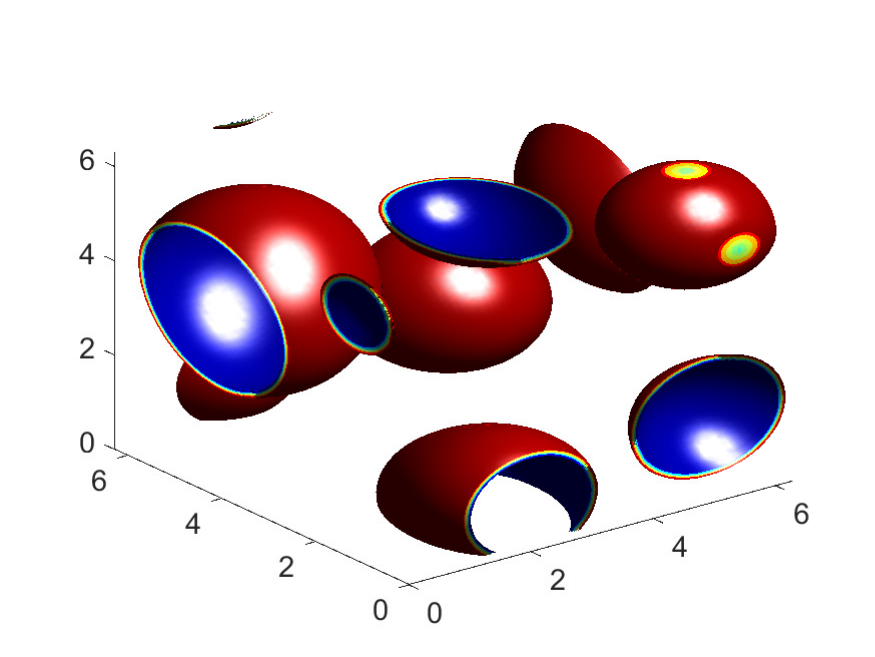}
\end{center}
\caption{
Numerical solution to the CHVM equation in 3D with small random initial data $(\bar u_0 =0.55)$ on $\Omega=[0,2\pi]^3$ using BHM-IMEX1 with $64^3$ elements and $h=0.01$. Simulation snapshots taken at $t=1,10,50,100$ (a-d respectively).  
}
\label{fig:CH3D}
\end{figure}

We also  consider the CH equation with variable mobility \eqref{CHVMeqn} with
$\cM(u)=1-\omega ^2u^2$,  $\cW(u)=1/4u^4-1/2u^2$ and $\omega \in [0,1]$. 
Here we choose $\omega=0.95$ as it was done in \cite{HectorC2013} which allows us to consider a case of variable mobility that is nearly degenerate \cite{DaiDu1,DaiDu2}.
We consider an initial condition in the form $u_0(x,y,z)=\bar u_0 + \eta u_1 $
where $\bar u_0$ is a constant with $-1<\bar u_0<1$,  and $u_1$ represents uniformly generated  random numbers between $(-1,1)$ and $\eta$ is a small parameter $(\eta \ll 1).$  We let the initial condition vary around $(\bar u_0(x,y,z) \approx 0.55)$ and consider the computational domain $\Omega=[0,2\pi]^3$ with $64^3$ modes. Using BHM-IMEX1 with $\Mb=0.5$, $\epsilon=0.1$ and $h=0.01$, we compute the numerical solution of the CHVM \eqref{CHVMeqn} in three-dimensions and report the results in Figure \ref{fig:CH3D}. Figure~\ref{fig:CH3D} shows snapshots for times $t=1,10,50,100$. The expected dynamics of the CHVM equation can be observed in Figure \ref{fig:CH3D} as the coarsening process takes place in which spheres (instead of circles for the 2-dimensional problem) evolve into configurations through mass-exchange leading to bigger spheres dominating over smaller counterparts. 
We also employed adaptive time step techniques \cite{Soderlind,Ganapa} to allow for high-accuracy computations during fast dynamics and efficiency during slow dynamics. A detailed discussion of these simulation results is presented in the Supplementary Material.

\section{Conclusions}

This paper presented new methods for the numerical solution of nonlinear fourth-order diffusion equations with variable mobility. The methods are based on a splitting approach (BHM) followed by different time-stepping discretizations and iterations. The main advantage of the splitting is that implicit terms are chosen to be linear constant coefficient operators in $u$ and computations become very efficient using Fourier pseudo-spectral methods. Our methods improve the accuracy of the biharmonic modified approach and achieve second order accuracy. We also provided a computational study of  the splitting parameters associated with the fourth-order and second-order splitting $\Mb$, and $\Mh$ which until now have received limited attention \cite{bertozzi,barrettblowey}. Their influence in terms of accuracy and energy-stability of solutions was presented and discussed using the TF equation with different test problems and the CH equation as it was presented in \cite{glasnerOrizaga}.
\par
The robustness of the methods makes them appropriate to handle thin film problems and variable mobility CH equations along with extensions of such models. The methods are very efficient and the computing requirements are comparable to the ones required by the original BHM method. We believe our approach could be a powerful tool, that is easy to implement, to study problems associated with nonlinear diffusion equations with variable mobility. We also believe that this paper has good potential for developing more fast and accurate methods for nonlinear higher-order phase field equations and to re-invigorate discussions and conversations in scientific communities interested in the TF and CH equations. 
\par
Further work will include a study of coarsening dynamics in 2D for the thin film equation, and in 3D for the Cahn-Hilliard equation with variable mobility. We plan to migrate our schemes to graphic processing unit (GPU) architecture \cite{LAM2019}, in particular for 3D simulations, to parallelize and efficiently perform the computations associated with large-scale dynamics. We are also interested in the phase field crystal (PFC) equation \cite{elder2004modeling,vignal2014energy,wise2009energy}, the block copolymer (BCP) equation \cite{ch-app1} and the Functionalized Cahn-Hilliard (FCH) equation \cite{promislowM,orizaga2023}. Problems to consider will be additional extensions or generalizations of Cahn-Hilliard equations \cite{miranville2017cahn} including coupled Cahn-Hilliard/Allen-Cahn systems, Cahn-Hilliard/thin film equations with non-conserved fluxes \cite{Thiele2018,Thiele2014}, and CH coupled to Navier-Stokes equations. These mentioned problems have applications in the bio-sciences and complex domains with curved surfaces \cite{Tobias2003} are possible, where fast computational methods are not available for the linear operator, which will give rise to new and interesting computational challenges.
A specific question of interest is how to select the stabilizing in $\Fi(u)$ through the $M_i$ coefficients to optimize the accuracy of computations. While we have focused on the dependence on $\Mb$ via \eqref{M1n}, careful selection of $\Mh$ as well as $M_0$ in \eqref{FimEqn} may yield valuable improvements. Some primarily work on tri-harmonic splitting, with an added $M_3\Delta^3 u$ term has already been applied to a sixth-order phase field model \cite{orizaga2023}.




\section*{Acknowledgment}
SO was supported through a Phillip Griffiths Research Assistant Professorship at Duke University and through start-up funds provided by New Mexico Tech. SO also acknowledges New Mexico IDeA Networks of Bio-medical Research Excellence (NM-INBRE) and the 2024 NM-INBRE NIRF Faculty Fellowship award. TW acknowledges support from grant NSF DMS 2008255. 

\section*{Declarations}
The authors declare that they have no known competing financial interests or personal relationships that could have appeared to influence the work reported in this paper. 

\section*{Author Contributions}
Both authors contributed equally in the preparation and completion of this manuscript. 

\section*{Data availability}
The codes used for these experiments are available on the freely accessible
GitHub repository
https://github.com/sauloorizaga/Thin-Film-Equations

\section*{Supplementary Material}
Supplementary data to this article can be found online at
the journal's website.

{\small 
\biboptions{sort&compress}

}





\end{document}